\documentclass[12pt]{amsart}

\usepackage{amsfonts}
\usepackage{amsmath}
\usepackage{amssymb}
\usepackage{mathrsfs}
\usepackage{color}
\usepackage{amsthm}
\usepackage[all]{xy}

\newtheorem{theorem}{Theorem}[section]
\newtheorem{proposition}[theorem]{Proposition}
\newtheorem{lemma}[theorem]{Lemma}
\newtheorem{corollary}[theorem]{Corollary}

\theoremstyle{definition}

\newtheorem{remark}[theorem]{Remark}
\newtheorem{case}{Case}

\newcommand{\map}[1]{\xrightarrow{#1}}

\newcommand{\abs}[1]{\lvert#1\rvert}    

 \newcommand{\Des}{\mathrm{Des}}
 \newcommand{\Peak}{\mathrm{Peak}}
 
 \newcommand{\ISh}{\mathrm{S}}
\newcommand{\comp}{\vDash}              
           
\newcommand{\sumsub}[1]{\sum_{\substack{#1}}} 

\newcommand{\ipartfrac}[1]{\lfloor \frac{#1}{2} \rfloor} 
\newcommand{\ipart}[1]{\lfloor #1/2 \rfloor} 
\newcommand{\ipartn}{\ipart{n}}  

 \newcommand{\peakint}{p_-}   
 \newcommand{\peakaug}{p_+}  
 \newcommand{\Peakint}{\Peak^0}  
 \newcommand{\Peakaug}{\Peak} 

\newcommand{\zetaS}{\zeta_\calS}

\newcommand{\calD}{\mathcal{D}}
\newcommand{\calL}{\mathcal{L}}
\newcommand{\calS}{\mathcal{S}}
\newcommand{\calQ}{\mathcal{Q}}

\newcommand{\QSym}{{\calQ}\Sym}

\newcommand{\SSym}{{\calS}\Sym}
\newcommand{\Sym}{{\mathit{Sym}}}
\def\H{{\mathcal H}}
\def\K{{\mathcal K}}

\def\field{\Bbbk}

\newcommand{\blue}[1]{\color{blue}#1}
\newcommand{\red}[1]{\color{red}#1}

\begin{document}
 \title[Canonical characters]{Canonical characters on quasi-symmetric
 functions and bivariate Catalan numbers}

 \author[M.~Aguiar and S.~Hsiao]{Marcelo Aguiar and Samuel K. Hsiao}

 \address{Department of Mathematics\\ Texas A\&M University\\
 College Station, TX 77843, USA}
 \email{maguiar@math.tamu.edu}
 \urladdr{http://www.math.tamu.edu/$\sim$maguiar}

\address{Department of Mathematics\\
 University of Michigan\\
  Ann Arbor, MI 48109, USA}
 \email{shsiao@umich.edu}
 \urladdr{http://www.math.lsa.umich.edu/$\sim$shsiao/}

 \date{July 23, 2004}
\thanks{Aguiar supported in part by NSF grant DMS-0302423}
 \thanks{Hsiao supported in part by the NSF Postdoctoral Research Fellowship}
 

 \begin{abstract} Every character on a graded connected Hopf algebra decomposes uniquely as a product of an even character and an odd character~\cite{ABS}.
 We obtain explicit formulas for the even and odd parts of the universal character on the Hopf algebra of quasi-symmetric functions. They can be described in terms of  Legendre's beta function evaluated at half-integers, or in terms of {\em bivariate Catalan numbers:}
 \[C(m,n)=\frac{(2m)!(2n)!}{m!(m+n)!n!}\,.\]
 Properties of characters and of quasi-symmetric functions are then used to derive several interesting identities among bivariate Catalan numbers and in particular
 among Catalan numbers and central binomial coefficients.
\end{abstract}

 \maketitle

\tableofcontents

\section{Introduction}\label{S:intro}

The numbers
\begin{equation}\label{E:Gnumbers}
C(m,n):=\frac{(2m)!(2n)!}{m!(m+n)!n!}=\frac{\binom{2m}{m}\binom{2n}{n}}{\binom{m+n}{n}}
\end{equation}
appeared in work of Catalan~\cite[p. 207]{Cat1},~\cite[Sections CV and CCXIV]{Cat2},~\cite[pp. 110-113]{Cat3}, von Szily~\cite[pp. 89--91]{VZ}, Riordan~\cite[Chapter 3, Exercise 9, p. 120]{Rio}, and recent work of Gessel~\cite{Ges2}. We call them {\em bivariate Catalan numbers}. They are integers (and except for $C(0,0)=1$, they are all even). Special cases include the
central binomial coefficients and the Catalan numbers:
\[C(0,n)=\binom{2n}{n} \text{ \ and \ }
\frac{1}{2}C(1,n)=\frac{1}{n+1}\binom{2n}{n}\,.\]
In turn, the bivariate Catalan numbers are special cases of the {\em super ballot numbers} of Gessel~\cite{Ges2}. 

The algebra $\QSym$ of quasi-symmetric functions was introduced in earlier work of Gessel~\cite{Ges} as a source of generating functions for Stanley's $P$-partitions~\cite{Sta72}; since then, the literature on the subject has become vast.   
The linear bases of $\QSym$ are indexed by {\em compositions} $\alpha$ of $n$. Two important bases are given by the {\em monomial} and {\em fundamental} quasi-symmetric functions $M_\alpha$ and $F_\alpha$; for more details, see~\cite{Ges},~\cite[Chapter 4]{Malv}, or~\cite[Section 7.19]{Sta99}.

In~\cite{ABS}, an important universal property of $\QSym$ was derived.
Consider the functional $\zeta:\QSym\to\field$ obtained by
specializing one variable of a quasi-symmetric function to $1$ and all 
other variables to $0$. On the monomial and  fundamental
bases of $\QSym$, this functional is given by

\begin{equation}\label{E:zeta-QSym}
\zeta(M_\alpha)\ =\ \zeta(F_\alpha)\ =\ \begin{cases}
       1 & \text{ if $\alpha=(n)$ or $(\,)$,}\\
       0& \text{otherwise.}\end{cases}
\end{equation}

A {\em character} on a  Hopf algebra $\H$ is a morphism of algebras $\varphi:\H\to\field$:
\[\varphi(ab)=\varphi(a)\varphi(b)\,,\quad \varphi(1)=1\,.\]
$\QSym$ is a graded connected Hopf algebra and the functional $\zeta$ is a character on $\QSym$. The universal property states that given any graded connected Hopf algebra $\H$ and a character
$\varphi:\H\to\field$, there exists a
 unique morphism of graded Hopf algebras $\Phi:\H\to\QSym$ 
making the following diagram commutative~\cite[Theorem 4.1]{ABS}:
\[\xymatrix{{\ \H\ }\ar[rr]^{\Phi}\ar[rd]_{\varphi} & & {\QSym}\ar[ld]^{\zeta}\\
& \field & }\]
For this reason, we refer to $\zeta$ as the {\em universal}
character on $\QSym$. There are other characters on $\QSym$
canonically associated to $\zeta$ that are of interest to us.
In spite of the simple definition of $\zeta$,
these characters encompass important combinatorial information.
Some of these were explicitly described and studied in~\cite{ABS},
and shown to be closely related to a Hopf subalgebra of $\QSym$
introduced by Stembridge~\cite{Ste97}, to the generalized 
Dehn-Sommerville relations, and to other combinatorial constructions. 
Other canonical characters, less easy to describe but of a more
fundamental nature, are the object of this paper.

We review other relevant background and constructions from~\cite{ABS}.

Let $\H$ be an arbitrary Hopf algebra. The {\em convolution product} of two linear functionals $\rho,\psi:\H\to\field$ is
\begin{equation*}\label{E:convolution}
\H\map{\Delta}\H\otimes
\H\map{\rho\otimes\psi}\field\otimes\field\map{m}\field\,,
\end{equation*}
where $\Delta$ is the coproduct of $\H$ and $m$ is the 
product of the base field. We denote the convolution
product by $\rho\psi$. The set of characters on any Hopf algebra is a group under the  convolution product. The unit element is $\epsilon:\H\to\field$, the counit map of $\H$. The inverse of a character $\varphi$
is $\varphi^{-1}:=\varphi\circ S$, where $S$ is the antipode of $\H$.

Suppose that $\H$ is graded, i.e., $\H=\oplus_{n\geq 0}\H_n$ and the structure maps of $\H$ preserve this decomposition. Then $\H$ carries a canonical
automorphism defined on homogeneous elements $h$ of degree $n$
by $h\mapsto \Bar{h}:=(-1)^n h$. If $\varphi$ is a functional on $\H$, we define
a functional $\Bar{\varphi}$ by $\Bar{\varphi}(h)=\varphi(\Bar{h})$.
The functional $\varphi$  is said to be {\em even} if
\[\Bar{\varphi}=\varphi\]
and it is said to be {\em odd} if it is invertible with respect to convolution and
\[\Bar{\varphi}=\varphi^{-1}\,.\]
Suppose now that $\H$ is graded and {\em connected}, i.e., $\H_0=\field\cdot 1$. One of the main results of~\cite{ABS} states that any character
$\varphi$ on $\H$ decomposes uniquely as a product of characters
\[\varphi=\varphi_+\varphi_-\]
with $\varphi_+$ even and $\varphi_-$ odd~\cite[Theorem 1.5]{ABS}.

The main purpose of this paper is to obtain explicit descriptions
for the canonical characters $\zeta_+$ and $\zeta_-$ of $\QSym$.
We find that the values of both characters are given in terms of bivariate Catalan
numbers (up to signs and powers of $2$). On the monomial basis,
the values are Catalan numbers and central binomial coefficients (Theorem~\ref{T:zetaM}).
On the fundamental basis, general bivariate Catalan numbers intervene (Theorem~\ref{T:zetaF}). The connection with Legendre's beta function is
given in Remark~\ref{R:integral}.
The proofs rely
on a number of identities for these numbers, of which some
are known and others are new. In turn, the general
properties of even and odd characters imply further identities
that these numbers must satisfy. We obtain in this way a large supply of identities for Catalan numbers and central binomial coefficients (Section~\ref{S:app1}) and for bivariate Catalan numbers (Sections~\ref{S:app2} and~\ref{S:permutations}). As one should expect, some of these identities  may also be obtained by more standard combinatorial arguments, at least once one is confronted with them. Our methods, however, yield the identities without any previous knowledge of their form. We mention here four of the most representative among the identities we derive:
\begin{equation}\tag{\ref{E:app-antipodeM}}
\sumsub{\alpha<\beta\\ a_1\text{ odd}}
 \frac{(-1)^{k_e(\alpha)}}{4^{\ipart{k_o(\alpha)}}} \binom{2\ipart{k_o(\alpha)}}{\ipart{k_o(\alpha)}}=0\,;
\end{equation}
\begin{equation}\tag{\ref{E:CG6}}
\sum_{j=0}^h C_3(j)C_3(h-j)  =2\sum_{j=0}^{h-1}C_2(j)C_1(h-1-j)\,;
\end{equation}
\begin{equation}\tag{\ref{E:allperms-}}
 \sum_{\sigma \in S_n} (-1)^{\peakint(\sigma)} C\bigl(\peakint(\sigma), \ipartn - \peakint(\sigma)\bigr) = 4^{\ipartn }\,;
\end{equation}
\begin{equation}\tag{\ref{E:catalan-gessel}}
2C(r)=\frac{1}{4^s}C(r,s+1)+\sum_{j=1}^s \frac{1}{4^j} C(r+1,j)\,.
\end{equation}

In~\eqref{E:app-antipodeM}, $k_e(\alpha)$ and $k_o(\alpha)$ are the number of even parts and the number of odd parts of a composition $\alpha$,  $\beta=(b_1,\ldots,b_k)$ is any fixed composition such  that $k_e(\beta)\equiv 0$ and $b_1\equiv b_k\mod 2$, and the sum is over  those compositions $\alpha=(a_1,\ldots,a_h)$ whose first part  is odd and which are  strictly refined by $\beta$.
The numbers $C_i(j)$ appearing in~\eqref{E:CG6}
are {\em central Catalan numbers}; see~\eqref{E:CGnumbers}. In~\eqref{E:allperms-}, $\peakint(\sigma)$ denotes the number of interior peaks of the permutation $\sigma$; see Section~\ref{S:permutations}. Equation~\eqref{E:catalan-gessel} expresses a Catalan number in terms of bivariate Catalan numbers.

In Section~\ref{S:inverses} we derive explicit formulas for the even and odd characters entering in the decomposition of the inverse  (with respect to convolution) of the universal character, and deduce some more identities, including~\eqref{E:catalan-gessel}.

We work over a field $\field$ of characteristic different from $2$.

\section{Even and odd characters}\label{S:evenodd}

Let $\H$ be a graded connected Hopf algebra and $\varphi:\H\to\field$
a linear functional such that $\varphi(1)=1$ (this holds  if $\varphi$ is a character). Let $\varphi_n$ denote the restriction of $\varphi$
to the homogeneous component of $\H$ of degree $n$. By assumption,
$\varphi_0=\epsilon_0$, where $\epsilon$ is the counit of $\H$. This guarantees that $\varphi$
is invertible with respect to convolution:
the inverse functional $\varphi^{-1}$ is determined
by the recursion
\[(\varphi^{-1})_n=-\sum_{i=1}^{n}\varphi_i(\varphi^{-1})_{n-i}\]
with initial condition $(\varphi^{-1})_0=\epsilon_0$.

\begin{lemma}\label{L:even-odd-decomposition}
Let $\H$ be a graded connected Hopf algebra and $\varphi:\H\to\field$ a linear
functional such that $\varphi(1)=1$. There are unique linear functionals $\rho,\psi:\H\to\field$ such that
\begin{itemize}
\item[(a)] $\rho(1)=\psi(1)=1$,
\item[(b)] $\rho$ is even and $\psi$ is odd,
\item[(c)] $\varphi=\rho\psi$.
\end{itemize}
Moreover, if $\varphi$ is a character then so are $\rho$ and $\psi$.
\end{lemma}
\begin{proof} Items (a), (b), and (c) can be derived as in the proof of~\cite[Theorem 1.5]{ABS}, while~\cite[Proposition 1.4]{ABS} guarantees that  if $\varphi$ is a character then so are $\rho$ and $\psi$.
\end{proof}
In this situation, we write $\varphi_+:=\rho$ and $\varphi_-:=\psi$ and  refer to them as the {\em even part} and the {\em odd part} of $\varphi$.
According to the results cited above, $\varphi_+$ is uniquely determined by the recursion
\[(-1)^n\varphi_n\ =\ 2(\varphi_+)_n+(\varphi^{-1})_n+
\sumsub{i+j+k=n\\0\leq i,j,k< n}(\varphi_+)_i(\varphi^{-1})_j(\varphi_+)_k\,,\]
and $\varphi_-$ by
\[(\varphi_-)_n=\varphi_n-\sum_{i=1}^{n}(\varphi_+)_i(\varphi_-)_{n-i}\]
with initial conditions $(\varphi_+)_0=(\varphi_-)_0=\epsilon$.

\begin{lemma}\label{L:even-odd-functorial}
Suppose $\H$ and $\K$ are graded connected Hopf algebras, $\varphi:\H\to\field$ and $\psi:\K\to\field$ are characters, and $\Phi:\H\to\K$ is a morphism of graded Hopf algebras such that
\[\xymatrix{
{\H}\ar[rr]^{\Phi}\ar[rd]_{\varphi} & & {\K}\ar[ld]^{\psi}\\
& {\field} & }\]
commutes. Then diagrams
\[\xymatrix{
{\H}\ar[rr]^{\Phi}\ar[rd]_{\varphi_+} & & {\K}\ar[ld]^{\psi_+}\\
& {\field} & } \qquad
\xymatrix{
{\H}\ar[rr]^{\Phi}\ar[rd]_{\varphi_-} & & {\K}\ar[ld]^{\psi_-}\\
& {\field} & }\]
commute as well.
\end{lemma}
\begin{proof} Composition with $\Phi$ gives a morphism from the character group of $\K$ to the character group of $\H$ which preserves the canonical involution $\varphi\mapsto\Bar{\varphi}$. Thus $\psi=\psi_+\psi_-$ implies
$\psi\circ\Phi=(\psi_+\circ\Phi)(\psi_-\circ\Phi)$,  $\psi_+\circ\Phi$ is even, and $\psi_-\circ\Phi$ is odd. By uniqueness in Lemma~\ref{L:even-odd-decomposition}, $\psi_+\circ\Phi=\varphi_+$ and $\psi_-\circ\Phi=\varphi_-$.
\end{proof}

When $\H=\QSym$ and $\varphi=\zeta$ is the universal character~\eqref{E:zeta-QSym}, we refer to $\zeta_+$ and $\zeta_-$ as the
{\em canonical} characters of $\QSym$.

Let  $\rho$ and $\psi$ be arbitrary characters on $\QSym$. 
For later use, we describe the convolution product $\rho\psi$ explicitly.
Given a composition  $\alpha=(a_1,\ldots,a_k)$ of a positive integer $n$
and $0\leq i\leq n$, let $\alpha_i = (a_1,\ldots, a_i)$ and $\alpha^i=(a_{i+1},\ldots,a_k)$. We agree that $\alpha_0=\alpha^k=(\,)$ (the empty composition). The coproduct of $\QSym$ is
\[\Delta(M_\alpha)=\sum_{i=0}^k M_{\alpha_i}\otimes M_{\alpha^i}\,,\]
where $M_{(\,)}$ denotes the unit element $1\in\QSym$. It follows that
\begin{equation}\label{E:convolution-QSym}
(\rho\psi)(M_\alpha)=\sum_{i=0}^k \rho(M_{\alpha_i})\psi(M_{\alpha^i})\,.
\end{equation}
The counit is
\[\epsilon(M_\alpha)=\begin{cases} 1 & \text{ if $\alpha=(\,)$,}\\
0 & \text{ otherwise.} \end{cases}\]
%
%

\section{The canonical  characters  of $\QSym$  on the
monomial basis} \label{S:can-even}

For any non-negative integer $m$, let 
\[B(m):= C(0,m)=\binom{2m}{m} \text{ \ and \ }
C(m) := \frac{1}{2}C(1,m)=\frac{1}{m+1}\binom{2m}{m}\,;\]
these are the central binomial coefficients and the Catalan numbers.

\begin{lemma}\label{L:classical-conv} For any non-negative integer $m$,
\begin{align} 
\label{E:classical-conv}
B(m) &= 2 \sum_{i=1}^m C(i-1) B(m-i)\,,\\
\label{E:classical-conv2}
2^{2m} &= \sum_{i=0}^m B(m) B(m-i)\,.
\end{align}
\end{lemma}
\begin{proof} These are well-known identities. They appear in~\cite[Formulas (3.90) and (3.92)]{Gou}, and~\cite[Chapter 3, Exercise 9, p. 120, and Section 4.2, Example 2, p. 130]{Rio}. For bijective proofs, see~\cite[Formulas (2) and (8)]{EK}.
\end{proof}

For a composition $\alpha$, let $\abs{\alpha}$ denote the sum of the parts of $\alpha$, $k(\alpha)$ the number of parts of $\alpha$, $k_e(\alpha)$  the number of even parts
of $\alpha$, and $k_o(\alpha)$ the number of odd parts of $\alpha$. Note that
\begin{equation}\label{E:oddparts}
k_o(\alpha)\equiv\abs{\alpha}\mod 2\,.
\end{equation}

\begin{theorem}\label{T:zetaM}
Let $\alpha=(a_1,\ldots,a_k)$ be a composition of a positive integer $n$. Then
\begin{eqnarray}
\label{E:zeta-M}
\zeta_{-}(M_\alpha) &=& \left\{\begin{array}{ll}
\displaystyle{ \frac{(-1)^{k_e(\alpha)}}{2^{2\ipart{k_o(\alpha)}}}
 C\left(0, \ipart{k_o(\alpha)} \right)}  & 
  \mbox{if $a_k$ is odd,}  \\
\rule{0pt}{20pt}0 & \mbox{if $a_k$ is even;}
\end{array} \right. 
\\
\label{E:zeta+M}
\zeta_{+}(M_\alpha)  &=&  \left\{\begin{array}{ll}
\displaystyle{\frac{(-1)^{k_e(\alpha)+1}}{2^{k_o(\alpha)}}
 C\left(1,k_o(\alpha)/2-1\right)} & \mbox{if $a_1$ and $a_k$ are odd and $n$ is even,} \\
 1 & \mbox{if $\alpha=(n)$ and $n$ is even,}\\
0 & \mbox{otherwise.}  \end{array}\right. 
\end{eqnarray}
We also have $\zeta_-(1) = \zeta_+(1) = 1.$
\end{theorem}
\begin{proof} Let $\rho,\psi:\QSym\to\field$ be the linear maps 
defined by the proposed formula for  $\zeta_{+}$ and
 $\zeta_{-}$, respectively. According to Lemma~\ref{L:even-odd-decomposition}, to conclude $\rho=\zeta_+$ and $\psi=\zeta_-$, it suffices to show that
$\rho\psi = \zeta$,  $\bar{\rho}=\rho$, and  $\bar{\psi}=\psi^{-1} $. 

Since $\rho$ vanishes on all compositions of $n$ when $n$ is odd, we have $\bar{\rho}=\rho$.

We  show that $\rho\psi = \zeta$.
Let $k_e:=k_e(\alpha)$ and $k_o:=k_o(\alpha)$.

\begin{case}
Suppose that $k=1$, so $\alpha = (n)$. We have
$\rho(M_{(n)}) = 1$ if $n$ is even and $0$ if $n$ is odd; also 
$\psi(M_{(n)}) = 0$ if $n$ is even and $1$ if $n$ is odd. Thus
$(\rho\psi)(M_{(n)}) = \rho(M_{(n)}) + \psi(M_{(n)}) = 1=\zeta(M_{(n)})$.
\end{case}

In all remaining cases $k>1$ and $\zeta(M_\alpha)=0$ by~\eqref{E:zeta-QSym}.

\begin{case}
Suppose that $k>1$ and $a_k$ is even. By~\eqref{E:convolution-QSym} we have
\[(\rho\psi)(M_\alpha)=\sum_{i=0}^{k-1}
\rho(M_{\alpha_i})\psi(M_{\alpha^i})+\rho(M_\alpha)=0\]
by the  second alternative of~\eqref{E:zeta-M} applied to $\alpha^i$ and the
 third alternative of~\eqref{E:zeta+M} applied to $\alpha$.
\end{case}

\begin{case}
Suppose that  $k>1$, $a_k$ is odd, and $a_1$ is even. In this case $\rho(M_{\alpha_i})=0$ for each $i>1$, so by~\eqref{E:convolution-QSym}
\begin{eqnarray*}
\lefteqn{ (\rho\psi)(M_\alpha) = 
\psi(M_\alpha) + \rho(M_{(a_1)}) \psi(M_{(a_2,\ldots,a_k)}) =} 
\\  & &
\frac{(-1)^{k_e}}{2^{2\lfloor k_o/2 \rfloor}}C(0,\lfloor k_o/2\rfloor) +
\frac{(-1)^{k_e-1}}{2^{2\lfloor k_o/2 \rfloor}}C(0,\lfloor k_o/2\rfloor)
=0.
 \end{eqnarray*}
 \end{case}
\begin{case} Suppose that $k>1$ and $a_k$, $a_1$, and $n$ are odd. By~\eqref{E:zeta+M}, we have
$\rho(M_{\alpha_i})=0$ unless $i=0$ or $a_i$ is odd and $\abs{\alpha_i}$ is even.
Hence
\begin{eqnarray*}
(-1)^{k_e}2^{k_o-1}(\rho\psi)(M_\alpha) &=& 
(-1)^{k_e}2^{k_o-1}\Bigl(\psi(M_\alpha) + \sum_{{1\le i\le k-1 \atop a_i \;\; \mathrm{ odd}}\atop |\alpha_i| \;\;\mathrm{even}  } 
  \rho(M_{\alpha_i})\psi(M_{\alpha^i})  \Bigr) \nonumber \\
&=& 
B\left(\frac{k_o(\alpha)-1}{2}\right)- \displaystyle{2 \sum_{{1\le i\le k-1 \atop a_i \;\; \mathrm{ odd}}\atop |\alpha_i| \;\;\mathrm{even} } 
C\left(\frac{k_o(\alpha_i)}{2}-1\right) B\left( \frac{k_o(\alpha^i)-1}{2} \right)\,.} 
\end{eqnarray*}
We used~\eqref{E:oddparts} and the fact that $\ipart{h}=\frac{h-1}{2}$ when $h$ is odd.

Deleting the even parts of $\alpha$ 
and changing every odd part of $\alpha$ to $1$ does not change the right-hand side of this equation. Thus we may assume without loss of generality that
$\alpha=(1,1,\ldots,1)\comp n = 2m+1.$ In this case,
\[
(-1)^{k_e}2^{k_o-1}(\rho\psi)(M_\alpha)=B(m) -2 \sum_{j=1}^m C(j-1) B(m-j)=0\,,
\]
by Lemma~\ref{L:classical-conv}.
\end{case}
\begin{case}\label{case5} Suppose that $k>1$, $a_k$ and $a_1$ are odd, and $n$ is even. Then
\begin{align*}
(-1)^{k_e}2^{k_o}(\rho\psi)(M_\alpha) &=  (-1)^{k_e}2^{k_o}\Bigl(\psi(M_\alpha) + \sumsub{1\le i\le k-1 \\ a_i \text{ odd} \\ \abs{\alpha_i} \text{ even}  }
  \rho(M_{\alpha_i})\psi(M_{\alpha^i})  + \rho(M_\alpha) \Bigr)  \\
 &=B\left(\frac{k_o(\alpha)}{2} \right)-2 \sumsub{1\le i\le k-1 \\ a_i \text{ odd}\\ \abs{\alpha_i} \text{ even} } 
C\left(\frac{k_o(\alpha_i)}{2}-1\right) B\left(\frac{k_o(\alpha^i)}{2}\right) -2C\left(\frac{k_o(\alpha)}{2} -1\right)\,. 
\end{align*} 
As before, it suffices to consider the special case
$\alpha=(1,1, \ldots,1) \comp 2m$. In this case,
\begin{align*}
(-1)^{k_e}2^{k_o}(\rho\psi)(M_\alpha) &=B(m) - 2 \sum_{j=1}^{m-1} C(j-1) B(m-j) -2C(m-1)\\
&=B(m) -2 \sum_{j=1}^m C(j-1) B(m-j)=0\,,
\end{align*}
again by Lemma~\ref{L:classical-conv}.
\end{case}
\medskip
It remains to show that $\psi^{-1} = \bar{\psi}$. Since $\psi$ is invertible ($\psi(1)=1$), it suffices to show that $\psi \bar{\psi} = \varepsilon$. 
\setcounter{case}{0}
\begin{case} If $a_k$ is even then $\bar{\psi}(M_{\alpha^i})=0$ for every $i<k$ and $\psi(M_{\alpha})=0$ by~\eqref{E:zeta-M}, so $(\psi\bar{\psi})(M_\alpha)=0$.
\end{case}
\begin{case} Suppose $a_k$  and $n$ are odd. In this case $\bar{\psi}(M_\alpha)=-\psi(M_{\alpha})$, so $(\psi\bar{\psi})(M_\alpha)=\sum_{i=1}^{k-1}\psi(M_{\alpha_i})\bar{\psi}(M_{\alpha^i})$. By~\eqref{E:oddparts}, $k_o(\alpha_i)+k_o(\alpha^i)=k_o(\alpha)$ is odd. Hence
\[\ipart{k_o(\alpha_i)}+\ipart{k_o(\alpha^i)}=\frac{k_o(\alpha)-1}{2}
\text{ \ and \ } k_e(\alpha_i)+k_e(\alpha^i)+\abs{\alpha^i}\equiv k(\alpha)+\abs{\alpha_i} \mod 2\,.\]
It follows that
\begin{eqnarray*}
(-1)^{k(\alpha)} 2^{k_o(\alpha)-1}(\psi\bar{\psi})(M_\alpha)
& = &  
 \sum_{1\le i \le k-1 \atop a_i \;\mathrm{odd}} (-1)^{|\alpha_i|} 
B(\ipart{k_o(\alpha_i)})B(\ipart{k_o(\alpha^i)})\,.
\end{eqnarray*}
As before, we may assume $\alpha = (1, 1,\ldots, 1) \comp n=2m+1$, in which case the above sum becomes
\[\sum_{i=1}^{2m}(-1)^i B(\ipart{i})B(\ipart{2m+1-i})\,.\]
As $j$ runs from $1$ to $m$, the terms in this sum corresponding to $i=2j$ and $i=2m-2j+1$ are $B(j)B(m-j)$ and $-B(m-j)B(j)$, respectively. Since this covers
all terms, this sum is zero.
\end{case}
\begin{case} Suppose $a_k$ is odd and $n$ is even.
Similar considerations lead to
\begin{eqnarray*}
\lefteqn{(-1)^{k(\alpha)} 2^{k_o(\alpha)}(\psi\bar{\psi})(M_\alpha) = } \\
 & &  B\Bigl(\frac{k_o(\alpha)}{2}\Bigr) + \sum_{{1\le i \le k-1 \atop a_i \;\mathrm{odd}}\atop |\alpha_i| \; \mathrm{even}} 
B\Bigl(\frac{k_o(\alpha_i)}{2}\Bigr)B\Bigl(\frac{k_o(\alpha^i)}{2}\Bigr)  - 
2^2\sum_{{1\le i \le k-1 \atop a_i \;\mathrm{odd}}\atop |\alpha_i| \; \mathrm{odd}} 
B\Bigl(\frac{k_o(\alpha_i)-1}{2}\Bigr)B\Bigl(\frac{k_o(\alpha^i)-1}{2}\Bigr).
\end{eqnarray*}
Once again, we may assume $\alpha = (1, 1,\ldots, 1) \comp n=2m$. Then showing that $(\psi\bar{\psi})(M_\alpha)=0$ is equivalent to showing that
$$\sum_{i=0}^{m} B(i) B(m-i) = 2^2 \sum_{i=0}^{m-1} B(i) B(m-1-i)\,.$$
This equality follows from~\eqref{E:classical-conv2} in Lemma~\ref{L:classical-conv}.
\end{case}
The proof is complete.\end{proof}

\section{Application: Identities for Catalan numbers and central binomial coefficients}\label{S:app1}

In the proof of Theorem~\ref{T:zetaM}, we did not need to show that the functionals defined by~\eqref{E:zeta-M}, \eqref{E:zeta+M} are characters (morphisms of algebras); indeed, this fact follows from our argument. We may derive interesting identities involving Catalan numbers or central binomial coefficients from this property. To this end, we first describe the product of two monomial quasi-symmetric functions. This result is known from~\cite[Lemma 3.3]{Eh96}, ~\cite{Ho00}, and~\cite[Section 5]{TU96}. We present here an equivalent but more convenient description due to Fares~\cite{Far}.

Given non-negative integers $p$ and $q$, consider the set $\calL(p,q)$ of lattice paths from $(0,0)$ to $(p,q)$ consisting of unit steps  which are either horizontal, vertical, or diagonal (usually called {\em Delannoy} paths). An element  of $\calL(p,q)$ is thus a sequence  $L=(\ell_1,\ldots,\ell_{s})$ such that each $\ell_i$ is either $(1,0)$, $(0,1)$, or $(1,1)$, and $\sum \ell_i=(p,q)$. Let $h$, $v$, and $d$ be the number of horizontal, vertical, and diagonal steps in $L$. Then $h+d=p$, $v+d=q$,
and $s=h+v+d=p+q-d$. The number of lattice paths in $\calL(p,q)$ with $d$
diagonal steps is the multinomial coefficient
\[\binom{s}{h,\,v,\,d}=\binom{p+q-d}{p-d,\,q-d,\,d}\,,\]
since such a path is determined by the decomposition of the set of steps into the subsets of horizontal, vertical, and diagonal steps.

Given compositions $\alpha=(a_1,\ldots,a_p)$ and $\beta=(b_1,\ldots,b_q)$ and $L\in\calL(p,q)$, we label each step of $L$ according to its horizontal and vertical projections, as indicated in the example below
($p=5$, $q=4$, and $L$ is in red):

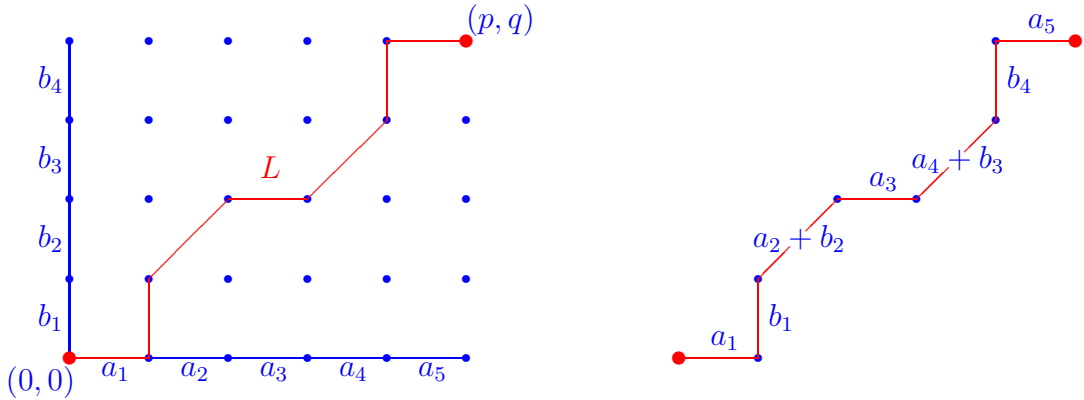
\[ \begin{picture}(160,140)(-15,-15){\color{blue}
   \put(0,0){\line(1,0){150}} \put(0,0){\line(0,1){120}}
   \put(-24,-12){$(0,0)$}\put(150,125){$(p,q)$}
    \put(12,-7){$a_1$} \put(42,-7){$a_2$} \put(72,-7){$a_3$} \put(102,-7){$a_4$} \put(132,-7){$a_5$}
    \put(-12,12){$b_1$} \put(-12,42){$b_2$} \put(-12,72){$b_3$} 
    \put(-12,102){$b_4$} 
    \put(0,30){\circle*{3}} \put(0,60){\circle*{3}} \put(0,90){\circle*{3}} \put(0,120){\circle*{3}} 
 \put(30,0){\circle*{3}}    \put(30,30){\circle*{3}} \put(30,60){\circle*{3}} \put(30,90){\circle*{3}} \put(30,120){\circle*{3}} 
 \put(60,0){\circle*{3}}    \put(60,30){\circle*{3}} \put(60,60){\circle*{3}} \put(60,90){\circle*{3}} \put(60,120){\circle*{3}} 
  \put(90,0){\circle*{3}}    \put(90,30){\circle*{3}} \put(90,60){\circle*{3}} \put(90,90){\circle*{3}} \put(90,120){\circle*{3}} 
 \put(120,0){\circle*{3}}    \put(120,30){\circle*{3}} \put(120,60){\circle*{3}} \put(120,90){\circle*{3}} \put(120,120){\circle*{3}} 
  \put(150,0){\circle*{3}}    \put(150,30){\circle*{3}} \put(150,60){\circle*{3}} \put(150,90){\circle*{3}} 
         }{\color{red}
     \put(0,0){\circle*{5}}   \put(150,120){\circle*{5}}   
      \put(0,0){\line(1,0){30}} \put(60,60){\line(1,0){30}} \put(120,120){\line(1,0){30}}
       \put(30,0){\line(0,1){30}}\put(120,90){\line(0,1){30}}
       \put(30,30){\line(1,1){30}} \put(90,60){\line(1,1){30}}
       \put(72,68){$L$}
         }\end{picture}
         \qquad\qquad\qquad
 \begin{picture}(160,140)(-15,-15){\color{blue}
  \put(12,5){$a_1$}  \put(34,12){$b_1$} 
  \put(28,42){$a_2+b_2$} \put(72,65){$a_3$} \put(88,72){$a_4+b_3$} 
  \put(124,102){$b_4$} \put(132,125){$a_5$}   
 \put(30,0){\circle*{3}}    \put(30,30){\circle*{3}} \put(60,60){\circle*{3}}
  \put(90,60){\circle*{3}} \put(120,90){\circle*{3}} \put(120,120){\circle*{3}} 
          }{\color{red}
     \put(0,0){\circle*{5}}   \put(150,120){\circle*{5}}   
      \put(0,0){\line(1,0){30}} \put(60,60){\line(1,0){30}} \put(120,120){\line(1,0){30}}
       \put(30,0){\line(0,1){30}}\put(120,90){\line(0,1){30}}
       \put(30,30){\line(1,1){12}} \put(90,60){\line(1,1){12}}
        \put(60,60){\line(-1,-1){12}} \put(120,90){\line(-1,-1){12}}
         }\end{picture}\]

Then we obtain a composition $q_L(\alpha,\beta)$ by reading off the labels
along the path $L$ in order. In the example above, 
\[q_L(\alpha,\beta)=(a_1,b_1,a_2+b_2,a_3,a_4+b_3,b_4,a_5)\,.\]
The composition $q_L(\alpha,\beta)$ is the {\em quasi-shuffle} of $\alpha$ and $\beta$ corresponding to $L$. If $L$ does not involve diagonal steps, then $q_L(\alpha,\beta)$ is an
ordinary shuffle.

The product of two monomial quasi-symmetric functions is given by
\begin{equation}\label{E:productQ}
M_\alpha\cdot M_\beta=\sum_{L\in\calL(p,q)} M_{q_L(\alpha,\beta)}\,.
\end{equation}

For our first application we make use of the fact that $\zeta_-$ is a character.

\begin{corollary}\label{C:central-prod} Let $n,m$ be non-negative integers  not both equal to $0$. Then
\begin{multline}\label{E:central-prod}
\sum_{d=0}^{\min(n,m)}\frac{(-1)^d}{4^{\ipart{(n+m-2d)}}}\frac{n+m-2d}{n+m-d}
\binom{n+m-d}{n-d,\,m-d,\,d}\binom{2\ipart{(n+m-2d)}}{\ipart{(n+m-2d)}}\\
=\frac{1}{4^{\ipart{n}+\ipart{m}}}\binom{2\ipart{n}}{\ipart{n}}\binom{2\ipart{m}}{\ipart{m}}\,.
\end{multline}
\end{corollary}
\begin{proof} Let $\alpha=(1,1,\ldots,1)\comp n$ and $\beta=(1,1,\ldots,1)\comp m$. The set of lattice paths $\calL(n,m)$ splits as
\[ \calL(n,m)=\calL_H(n,m)\sqcup \calL_V(n,m)\sqcup \calL_D(n,m)\]
 according to whether the last step of the path is horizontal, vertical, or diagonal.

Choose $L\in \calL(n,m)$, let $d$ be the number of diagonal steps of $L$, and
$\gamma:=q_L(\alpha,\beta)$.

If $L\in\calL_D(n,m)$, the last part of $\gamma$ is $2$, so
$\zeta_-(M_\gamma)=0$ by~\eqref{E:zeta-M}. On the other hand, if $L\in\calL_H(n,m)\sqcup \calL_V(n,m)$ the last part of $\gamma$ is $1$, $k_e(\gamma)=d$ ($d$ parts are equal to $2$), and
$k_o(\gamma)=n+m-2d$ (all other parts are equal to $1$). Hence, by~\eqref{E:zeta-M},
\[\zeta_-(M_\gamma)= \frac{(-1)^{d}}{2^{2\ipart{(n+m-2d)} }}
 C(0, \ipart{(n+m-2d)})\]
Now, the number of paths in $\calL_H(n,m)\sqcup \calL_V(n,m)$ with $d$ diagonal steps is
\[\binom{n-1+m-d}{n-1-d,\,m-d,\,d}+\binom{n+m-1-d}{n-d,\,m-1-d,\,d}=\frac{n+m-2d}{n+m-d}\,\binom{n+m-d}{n-d,\,m-d,\,d}\,.\]
Therefore, by~\eqref{E:productQ},
\[\zeta_-(M_\alpha\cdot M_\beta)=\sum_{d=0}^{\min(n,m)}\frac{(-1)^{d}}{2^{2\ipart{(n+m-2d)} }}\frac{n+m-2d}{n+m-d}\,\binom{n+m-d}{n-d,\,m-d,\,d} C(0, \ipart{(n+m-2d)})\,.\]
On the other hand, 
\[\zeta_-(M_\alpha)\cdot \zeta_-(M_\beta)=\frac{1}{2^{2\ipart{n} }}
 C(0, \ipart{n})\cdot \frac{1}{2^{2\ipart{m} }}
 C(0, \ipart{m})\,.\]
Since $\zeta_-$ is a character, the last two quantities are equal. Equating them
results in~\eqref{E:central-prod}.
\end{proof}

The special case of~\eqref{E:central-prod} when $n\geq m=1$ is
\begin{equation*}\label{E:central-prod1}
\frac{n+1}{4^{\ipart{(n+1)}}}\binom{2\ipart{(n+1)}}{\ipart{(n+1)}}
-\frac{n-1}{4^{\ipart{(n-1)}}}\binom{2\ipart{(n-1)}}{\ipart{(n-1)}}=
\frac{1}{4^{\ipart{n}}}\binom{2\ipart{n}}{\ipart{n}} 
\end{equation*}
which may be easily verified. On the other hand, when $n=m\geq 1$, ~\eqref{E:central-prod} becomes
\begin{equation*}\label{E:central-prod2}
\sum_{d=0}^{n-1}\frac{(-1)^d}{4^{n-d}}
\binom{2n-d-1}{d}\binom{2(n-d)}{n-d}^2\\
=\frac{1}{4^{2\ipart{n}}}\binom{2\ipart{n}}{\ipart{n}}^2\,.
\end{equation*}

\medskip

We now make use of the fact that $\zeta_+$ is a character to derive an identity involving Catalan numbers.

\begin{corollary}\label{C:catalan-prod} Let $n,m$ be positive integers not both equal to $1$ and such that $n\equiv m\mod 2$. Then
\begin{multline}\label{E:catalan-prod}
\sum_{d=0}^{\min(n,m)}(-1)^{d+1}2^{2d-1}\frac{(n+m-2d)}{(n+m-d)}\,\frac{(n+m-2d-1)}{(n+m-d-1)}\binom{n+m-d}{n-d,\,m-d,\,d}C\bigl((n+m)/2-d-1\bigr)\\
=\begin{cases} C(n/2-1) C(m/2-1) & \text{ if $n$ and $m$ are even,}\\
0 &  \text{ if $n$ and $m$ are odd.}
\end{cases}
\end{multline}
\end{corollary}
\begin{proof} Let $\alpha=(1,1,\ldots,1)\comp n$ and $\beta=(1,1,\ldots,1)\comp m$.  As in the proof of Corollary~\ref{C:central-prod}, we analyze which lattice
paths $L\in\calL(n,m)$ contribute to $\zeta_+(M_\alpha\cdot M_\beta)$. 
Note that since $(n,m)\neq (1,1)$, the second alternative in~\eqref{E:zeta+M}
never occurs. A path which starts or ends with a diagonal step does not contribute, by the third alternative in~\eqref{E:zeta+M}. For all  remaining paths $L$, if $d$ is the number of diagonal steps, then
\[\zeta_+(M_{q_L(\alpha,\beta)})=\frac{(-1)^{d+1}}{2^{n+m-2d-1}}
 C\bigl((n+m-2d)/2-1\bigr)\]
by the first alternative in~\eqref{E:zeta+M}. These paths start and end with a horizontal or vertical step, so their number  is
\begin{multline*}
2\binom{n-1+m-1-d}{n-1-d,\,m-1-d,\,d}+\binom{n-2+m-d}{n-2-d,\,m-d,\,d}+\binom{n+m-2-d}{n-d,\,m-2-d,\,d}\\
=\frac{(n+m-2d)}{(n+m-d)}\,\frac{(n+m-2d-1)}{(n+m-d-1)}\,\binom{n+m-d}{n-d,\,m-d,\,d}\,.
\end{multline*}
Applying $\zeta_+$ to both sides of~\eqref{E:productQ} leads to~\eqref{E:catalan-prod}.
\end{proof}

Suppose $m=1$, $n=2k+1$, $k\geq 1$. In this case~\eqref{E:catalan-prod} boils down to the simple identity
\[C(k)=\frac{2(2k-1)}{k+1}C(k-1)\,.\]
If $n=m$, the last term in the sum~\eqref{E:catalan-prod} is $0$, because of the factor
$n+m-2d$ (so there is no need to evaluate the Catalan number in this case). The formula becomes
\[\sum_{d=0}^{n-1}(-1)^{d+1}4^d(2n-2d-1)\binom{2n-d-2}{d}C(n-d-1)^2
=\begin{cases} C(n/2-1) ^2 & \text{ if $n$ is even,}\\
0 &  \text{ if $n$ is odd, $n>1$.}
\end{cases}\]

\medskip

In the proof of Theorem~\ref{T:zetaM} we established that formula~\eqref{E:zeta-M} defines an odd character by showing that $\zeta_-\Bar{\zeta}_-=\epsilon$.
Rewriting this property in terms of  the antipode $S$ of $\QSym$ leads to
new combinatorial identities that we analyze next. First, recall that $S$ is given by
\begin{equation}\label{E:antipodeM}
S(M_{\beta})=(-1)^{k(\beta)}\sum_{\alpha\leq\tilde{\beta}}M_\alpha\,,
\end{equation}
where $\tilde{\beta}=(b_k,\ldots,b_2,b_1)$ is the reversal of $\beta=(b_1,b_2,\ldots,b_k)$ and $\alpha\leq\gamma$ indicates that $\gamma$ is a refinement of $\alpha$~\cite[Proposition 3.4]{Eh96}; ~\cite[Corollaire 4.20]{Malv}.

\begin{corollary}\label{C:app-antipodeM} For any composition $\beta=(b_1,\ldots,b_k)$,
\begin{equation}\label{E:preapp-antipodeM}
\sumsub{\alpha\leq\beta\\ a_1\text{ odd}}
 \frac{(-1)^{k_e(\alpha)}}{4^{\ipart{k_o(\alpha)}}} \binom{2\ipart{k_o(\alpha)}}{\ipart{k_o(\alpha)}}=
  \begin{cases}
{\displaystyle \frac{1}{4^{\ipart{k_o(\beta)}}} \binom{2\ipart{k_o(\beta)}}{\ipart{k_o(\beta)}} }& \text{ if $b_k$ is odd,} \\
     0    & \text{ if $b_k$ is even;}
\end{cases}
\end{equation}
the sum being over those compositions $\alpha=(a_1,\ldots,a_h)$ whose first part  is odd and which are  refined by $\beta$.
\end{corollary}
\begin{proof} Since $\zeta_-$ is an odd character, we have
\[\zeta_-\circ S=(\zeta_-)^{-1}=\Bar{\zeta}_-\,.\]
Therefore, for any composition $\beta$ of $n$ we have, by~\eqref{E:antipodeM},
\[(-1)^{k(\beta)}\sum_{\alpha\leq\tilde{\beta}}\zeta_-(M_\alpha)=(-1)^n\zeta_-(M_\beta)\,.\]
Hence,  by~\eqref{E:zeta-M},
\[(-1)^{k(\beta)}\sumsub{\alpha\leq\tilde{\beta}\\ a_h\text{ odd}}
 \frac{(-1)^{k_e(\alpha)}}{2^{2\ipart{k_o(\alpha)} }} C(0, \ipart{k_o(\alpha)})=
  \begin{cases}
{\displaystyle \frac{(-1)^{n+k_e(\beta)}}{2^{2\ipart{k_o(\beta)}}} C(0,{\ipart{k_o(\beta)}}) }& \text{ if $b_k$ is odd,} \\
     0    & \text{ if $b_k$ is even.}
\end{cases}\]
By~\eqref{E:oddparts}, $k(\beta)=k_o(\beta)+k_e(\beta)\equiv n+k_e(\beta)$.
Replacing $\alpha$ for $\tilde{\alpha}$ in the left-hand side gives~\eqref{E:preapp-antipodeM}.
\end{proof}
Suppose that $k_e(\beta)\equiv 0 \mod 2$ and $b_1\equiv b_k\mod 2$.
In this case,  the summand in~\eqref{E:app-antipodeM} corresponding to $\alpha=\beta$ cancels with the right-hand side. We deduce that for any such composition $\beta$,
\begin{equation}\label{E:app-antipodeM}
\sumsub{\alpha<\beta\\ a_1\text{ odd}}
 \frac{(-1)^{k_e(\alpha)}}{4^{\ipart{k_o(\alpha)}}} \binom{2\ipart{k_o(\alpha)}}{\ipart{k_o(\alpha)}}=0\,.
 \end{equation}

Suppose now that $\beta=(1,1,\ldots,1)$. Consider the class of compositions of $n$ in which $r$ parts are
 odd (including the first part) and $s$ parts are even. Such a composition $\alpha$ is
 determined by the choice of a composition $(c_1,\ldots,c_{r+s})$ of $\frac{n+r}{2}$, plus the choice of $r-1$ parts out of the parts $c_1,\ldots,c_{r+s-1}$. In fact, once these choices are made, we obtain $\alpha$ by replacing
$c_{1}$ and the $r-1$ chosen $c_i$'s by $2c_i-1$, and replacing the
$s$ remaining $c_i$'s by $2c_i$. Any such $\alpha$ arises in this way uniquely.
Therefore, the number of such compositions is
\[\binom{\frac{n+r}{2}-1}{r+s-1}\binom{r+s-1}{r-1}\,.\]
Combining this with~\eqref{E:app-antipodeM} we deduce
\[\sum_{(r,s)\in T_n}\frac{(-1)^s}{4^{\ipart{r}}}\binom{\frac{n+r}{2}-1}{r+s-1}\binom{r+s-1}{r-1}\binom{2\ipart{r}}{\ipart{r}}=0\,;\]
the sum being over the set 
\[T_n:=\{(r,s) \mid r\equiv n\!\! \mod 2,\ r+2s\leq n,\ 0<r<n,\ 0\leq s\}\,.\]

We remark that this identity is less interesting than it may seem, since for each  $r<n$, ${\displaystyle \sum_{s=0}^{\frac{n-r}{2}}(-1)^s\, \binom{\frac{n+r}{2}-1}{r+s-1}\binom{r+s-1}{r-1}=0}$ by a Vandermonde convolution.

%
%

\section{The canonical  characters of $\QSym$  on the fundamental
basis}\label{S:can-odd}

For a composition $\alpha = (a_1, \ldots, a_k)$ of a non-negative integer $n$,
define
\begin{eqnarray}
\label{E:p-}
\peakint(\alpha)& = & \#\{i \ne k: a_i >1\}, \\
\label{E:p+}
\peakaug(\alpha) & = & \left\{\begin{array}{ll} 1+ \#\{i \ne 1,k : a_i >1\} & \mbox{if $k>1$,} \\
0 & \mbox{if $k\le 1$.}\end{array}\right.
\end{eqnarray}
A combinatorial interpretation for these two statistics in terms of peaks of permutations is given in Section~\ref{S:permutations}.

A composition $\alpha = (a_1, \ldots, a_k)$ may be conveniently represented by a {\em ribbon diagram}: a sequence of rows of squares, each row consisting of
 $a_i$ squares, and with the first square in row $i+1$ directly below the last
 square in row $i$. For instance the diagram

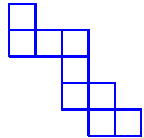
\[ \begin{picture}(60,50){\color{blue}
   \put(10,50){\line(1,0){10}} \put(10,40){\line(1,0){30}}
    \put(10,30){\line(1,0){30}} \put(30,20){\line(1,0){20}}
     \put(30,10){\line(1,0){30}} \put(40,0){\line(1,0){20}}
     \put(10,30){\line(0,1){20}}
    \put(20,30){\line(0,1){20}} \put(30,10){\line(0,1){30}}
    \put(40,0){\line(0,1){40}} \put(50,0){\line(0,1){20}}
    \put(60,0){\line(0,1){10}}
    }\end{picture}\]
 represents the composition $(1,3,1,2,2)$.
Note that $\peakint(\alpha)$ is the number of upper corners in the ribbon diagram of $\alpha$. To get a similar interpretation for $\peakaug(\alpha)$ one may augment the ribbon diagram of $\alpha$ by drawing an extra square to the left of the first row. Then $\peakaug(\alpha)$ is the number of upper corners in the augmented diagram. For instance, $\peakint(1,3,1,2,2)=2$ and $\peakaug(1,3,1,2,2)=3$, as illustrated below.
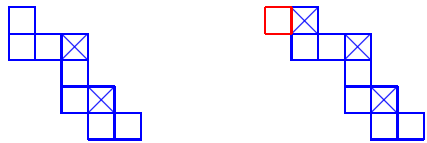
\[ \begin{picture}(60,50){\color{blue}
   \put(10,50){\line(1,0){10}} \put(10,40){\line(1,0){30}}
    \put(10,30){\line(1,0){30}} \put(30,20){\line(1,0){20}}
     \put(30,10){\line(1,0){30}} \put(40,0){\line(1,0){20}}
     \put(10,30){\line(0,1){20}}
    \put(20,30){\line(0,1){20}} \put(30,10){\line(0,1){30}}
    \put(40,0){\line(0,1){40}} \put(50,0){\line(0,1){20}}
    \put(60,0){\line(0,1){10}}
    \put(30,30){\line(1,1){10}}\put(40,30){\line(-1,1){10}}
    \put(40,10){\line(1,1){10}}\put(50,10){\line(-1,1){10}}
    }\end{picture}
    \qquad\qquad
    \begin{picture}(60,50){\color{blue}
   \put(10,50){\line(1,0){10}} \put(10,40){\line(1,0){30}}
    \put(10,30){\line(1,0){30}} \put(30,20){\line(1,0){20}}
     \put(30,10){\line(1,0){30}} \put(40,0){\line(1,0){20}}
    \put(10,30){\line(0,1){10}}
    \put(20,30){\line(0,1){20}} \put(30,10){\line(0,1){30}}
    \put(40,0){\line(0,1){40}} \put(50,0){\line(0,1){20}}
    \put(60,0){\line(0,1){10}}
    \put(10,40){\line(1,1){10}}\put(20,40){\line(-1,1){10}}
    \put(30,30){\line(1,1){10}}\put(40,30){\line(-1,1){10}}
    \put(40,10){\line(1,1){10}}\put(50,10){\line(-1,1){10}}
    }{\color{red}
     \put(0,40){\line(0,1){10}} \put(0,50){\line(1,0){10}} \put(0,40){\line(1,0){10}}
      \put(10,40){\line(0,1){10}}
     }\end{picture}\]

\medskip

The fundamental and monomial bases of $\QSym$ are related by
\begin{equation}\label{E:FM}
M_\alpha = \sum_{\beta\ge\alpha} (-1)^{k(\beta)-k(\alpha)} F_\beta\,,
\end{equation}
where, as before, for compositions  $\alpha$ and $\beta$ of $n$, $\beta\geq\alpha$ indicates that $\beta$ is refinement of $\alpha$.

\begin{theorem} \label{T:zetaF}
Let $\alpha$ be a composition of a non-negative integer $n$. Then
\begin{eqnarray}
\zeta_{-}(F_\alpha) & = & 
\displaystyle{\frac{(-1)^{\peakint(\alpha)}}{2^{2\ipartn}} 
C\left(\peakint(\alpha),\lfloor n/2 \rfloor-\peakint(\alpha)\right);} 
\label{E:zeta-F}\\
\zeta_{+}(F_\alpha) & = &  
\left\{\begin{array}{ll}
\displaystyle{\frac{(-1)^{\peakaug(\alpha)}}{2^{n}} 
C\left(\peakaug(\alpha),n/2 -\peakaug(\alpha)\right)}
& \mbox{ if $n$ is even,} \\
0 & \mbox{ if $n$ is odd.} \end{array} \right.
\label{E:zeta+F}
\end{eqnarray}
\end{theorem}

\begin{remark}\label{R:integral} Up to a sign and a factor of $\pi$, the above are values of Legendre's beta function at half-integers. Indeed, as already remarked by Catalan~\cite[Section CV]{Cat2},
\begin{equation}\label{E:integral}
\frac{1}{2^{2(p+q)}}C(p,q)=\frac{1}{\pi}B\Bigl(p+\frac{1}{2},q+\frac{1}{2}\Bigr)=\frac{1}{\pi}\int_0^{\pi}\sin^{2p}\theta\cos^{2q}\theta\,d\theta=\frac{1}{\pi}\int_0^1\frac{(1-u)^p u^q}{\sqrt{(1-u)u}}\,du\,.
\end{equation}
\end{remark}
\begin{remark}\label{R:power2}
Catalan shows~\cite[Section CCXIV]{Cat2},~\cite[pp. 110-113]{Cat3} that $C(p,q)$ is an integer, and that the highest power of $2$ that divides it equals the number of $1$'s in the binary decomposition of $p+q$. Equivalently,
the above fraction may be reduced as follows
\begin{equation}\label{E:power2}
\frac{1}{2^{2(p+q)}}C(p,q)=\frac{N}{2^k} \text{ \ with $N$ odd and }
k=\sum_{i\geq 0}\lfloor\frac{p+q}{2^i}\rfloor\,.
\end{equation}
\end{remark}
\smallskip
The proof of Theorem~\ref{T:zetaF} is given below. We  need two more statistics. For $\alpha$ as above, 
 let  
 \[u(\alpha) = \#\{ i \ne 1 : a_i >1\} \text{ \  and \ }
v(\alpha) = \#\{ i : a_i > 1\}\,.\] 

\begin{lemma} \label{lem:signs}
Let $m, j$ be non-negative integers. Then
\begin{itemize}
\item[(a)] $\displaystyle{\sum_{\gamma \comp m\atop v(\gamma) = j} (-1)^{k(\gamma)} 
= (-1)^{m+j} \binom{\lfloor m/2\rfloor}{j}}$;
\item[(b)] $\displaystyle{\sum_{\gamma\comp m \atop u(\gamma) = j} (-1)^{k(\gamma)}
= \left\{ \begin{array}{ll}  0 & \mbox{if $m$ is even,} \\
\displaystyle{(-1)^{m+j} \binom{\lfloor m/2\rfloor}{j}} & \mbox{if $m$ is odd.}
\end{array} \right.
}$
\end{itemize}
\end{lemma}
We thank Ira Gessel for supplying this  proof.
\begin{proof}
The generating
function for the left-hand side of (a) multiplied by $y^jx^m$ is just the
sum of  weights of all
compositions, where a part $a>1$ is weighted by $-yx^a$ and
a part equal to $1$ is weighted by $-x.$ This sum is
 $$\frac{1}{1- ( -x -yx^2-yx^3+\cdots)} 
 Ê=\frac{1}{1+x+\frac{yx^2}{1-x}}=\frac{1-x}{1-(1-y)x^2}.
$$
By expanding in powers of $x$ and $y$ we obtain the identity (a).
 For (b), the first part $a$ Êis weighted $-x^a$ no matter what
 $a$ is, so the generating function is 
 $-x/(1-x)$ times the generating function for (a), which is 
 $-x/ (1-(1-y)x^2),$ which  is the odd powers of $x$ in (a). 
\end{proof}

\begin{lemma} \label{L:G-convolve}
Let $i,j,m$ be non-negative integers. Then 
$$2^{2m} C(i,j) = \sum_{b=0}^m \binom{m}{b} C(i+b, m+j-b).$$
\end{lemma}
\begin{proof} In view of~\eqref{E:integral}, the proposed equality is equivalent to
\[\int_0^{\pi}\sin^{2i}\theta\cos^{2j}\theta\,d\theta=
 \sum_{b=0}^m \binom{m}{b}\int_0^{\pi}\sin^{2(i+b)}\theta\cos^{2(m+j-b)}\theta\,d\theta\,.\]
 This holds since
 \[ \sum_{b=0}^m \binom{m}{b}\sin^{2b}\theta\cos^{2(m-b)}\theta=(\sin^2\theta+\cos^2\theta)^m=1\,.\]
\end{proof}

To facilitate the proof of Theorem~\ref{T:zetaF}, we define
\begin{eqnarray*}
H_-(\alpha) 
& = & 
\sum_{\beta\ge\alpha}  (-1)^{k(\beta) + \peakint(\beta)+1}
C\left(\peakint(\beta),\ipart{n}-\peakint(\beta)\right) \;\;\;\;\; \mbox{for $\alpha\comp n$;} 
\\
H_+(\alpha) 
& = & 
\sum_{\beta\ge\alpha}  (-1)^{k(\beta) + \peakaug(\beta)+1}
C\left(\peakaug(\beta),n/2-\peakaug(\beta)\right) \;\;\; \mbox{for $\alpha\comp n$, $n$ even.}
\end{eqnarray*}


\begin{lemma}\label{L:H-alpha}
Suppose that $\alpha = (a_1,\ldots, a_k) \comp n$.
Then
$$
H_-(\alpha) = \left\{\begin{array}{ll}
(-1)^{n-1} 2^{n- k_o(\alpha) } C(0,\ipart{k_o(\alpha)}) & \mbox{if $a_k$ is odd,} \\
0 & \mbox{if $a_k$ is even.}
\end{array}\right.
$$
\end{lemma}
\begin{proof}
We will consider the cases $\alpha = (n)$ and $\alpha \ne (n)$ separately.

\setcounter{case}{0}
\begin{case} Suppose that $\alpha = (n)$. We need to show that
$H_-((n)) = 2^{n-1}$ if $n$ is odd and $H_-((n)) = 0$ if $n$ is even.
For a composition $\gamma$ and a positive integer $i$, let $\gamma i$ denote 
the concatenation of $\gamma$ with $(i)$. Using the fact that $\peakint(\gamma i) = 
v(\gamma),$ we have
\begin{eqnarray*}
H_-((n)) & = & \sum_{\beta \comp n} (-1)^{k(\beta) + \peakint(\beta)+1} 
                 C(\peakint(\beta),\ipartn -\peakint(\beta)) \\
                 & = & C(0,\ipartn) + \sum_{i=1}^{n-1} \sum_{\gamma \comp n-i} 
                 (-1)^{k(\gamma i)+ \peakaug(\gamma i)+1} 
                 C(\peakint(\gamma i), \ipartn - \peakint(\gamma i)) \\
                 & = & C(0,\ipartn) + \sum_{i=1}^{n-1} \sum_{\gamma \comp n-i}
                 (-1)^{k(\gamma) + v(\gamma)}
                 C(v(\gamma), \ipartn - v(\gamma)) \\
                 & = & C(0,\ipartn) + \sum_{i=1}^{n-1} \sum_{b=0}^{\ipart{(n-i)}}
                 \sum_{\gamma \comp n-i \atop v(\gamma) = b} 
                 (-1)^{k(\gamma) + b} C(b, \ipartn - b) \\
                 & = &  C(0,\ipartn) + \sum_{i=1}^{n-1} \sum_{b=0}^{\ipart{(n-i)}}
                 (-1)^b C(b, \ipartn -b) \sum_{\gamma\comp n-i \atop v(\gamma) = b}
                 (-1)^{k(\gamma)}.
\end{eqnarray*}
Make the substitution $m=n-i$ and apply
Lemma~\ref{lem:signs}(a) to get
\begin{eqnarray*}
H_-((n)) & = & C(0,\ipartn) + \sum_{i=1}^{n-1} \sum_{b=0}^{\ipart{m}}
                 (-1)^m \binom{\ipart{m}}{b} C(b, \ipartn-b)  \\
                 & = & C(0,\ipartn) + \sum_{b=0}^{\ipart{n-1}} C(b,\ipartn-b) 
                 \sum_{m=1}^{n-1} (-1)^m \binom{\ipart{m}}{b}.                 
\end{eqnarray*}
Let $S_b=\sum_{m=1}^{n-1} (-1)^m \binom{\ipart{m}}{b}$. Note that
when $m$ is even, $(-1)^m \binom{\ipart{m}}{b} + (-1)^{m+1}\binom{\ipart{(m+1)}}{b} = 0.$
So all but at most one term in the sum $S_b$ will cancel.

If $n$ is even, then $S_b = 0$ when $b\ge 1$ and $S_0 = -1$; hence
$H_-((n)) = C(0,\ipartn) + [-C(0,\ipartn)] = 0.$
If $n$ is odd,  
then $S_b = \binom{\ipartn}{b}$ when $b\ge 1$ and $S_0 = 0$; hence
by Lemma~\ref{L:G-convolve} (with $i=j=0$),
\begin{eqnarray*}
H_-((n)) &=& C(0,\ipartn) + \sum_{b=1}^{\ipartn} C(b,\ipartn - b) \binom{\ipartn}{b} \\
&=& \sum_{b=0}^{\ipartn} \binom{\ipartn}{b} C(b,\ipartn -b) \\
& = & 2^{2\ipart{n}}.
\end{eqnarray*}
Moreover, we have $2\ipartn  = n-1$ since $n$ is odd.
\end{case}

\begin{case}Suppose that $\alpha \ne (n)$.
If $\beta$ and $\gamma$ are 
non-empty compositions and $\beta\gamma$ denotes their concatenation,
then $\peakint(\beta\gamma) = v(\beta) + u(\tilde{\gamma})$, 
where $\tilde{\gamma}$ is the reversal of $\gamma$. 
Using this observation we may write
\begin{eqnarray}
H_-(\alpha) &=& 
\sum_{\beta\ge \alpha_{k-1} \atop \gamma \comp a_k} 
(-1)^{k(\beta\gamma) + \peakint(\beta\gamma)} 
C(\peakint(\beta\gamma), \ipartn - \peakint(\beta\gamma))  
\nonumber \\
& = &
\sum_{\beta \ge \alpha_{k-1}} 
(-1)^{k(\beta)+v(\beta)} \sum_{\gamma \comp a_k} (-1)^{k(\gamma)+u(\tilde{\gamma})} 
C(v(\beta) + u(\tilde{\gamma}),\ipartn-v(\beta) - u(\tilde{\gamma}))
\nonumber \\
& = &
\sum_{\beta \ge \alpha_{k-1}} (-1)^{k(\beta)+v(\beta)}
 \sum_{b=0}^{\lfloor a_k/2\rfloor} (-1)^{b} C(v(\beta)+b, \ipartn - v(\beta)-b) 
\sum_{\gamma\comp a_k \atop u(\tilde{\gamma})=b} (-1)^{k(\gamma)} \,.
\label{E:H-1} 
\end{eqnarray}
By Lemma~\ref{lem:signs}(b), the second sum in  \eqref{E:H-1} equals
\begin{eqnarray} 
\lefteqn{ 
\left\{\begin{array}{ll}
\displaystyle{(-1)^{a_k} \sum_{b=0}^{\lfloor a_k/2\rfloor}  \binom{\lfloor a_k/2\rfloor}{b}
C\left(v(\beta) + b, \ipartn - v(\beta)-b\right)}
    & \mbox{if $a_k$ is odd,} \\
0 
    & \mbox{if $a_k$ is even}
\end{array}\right. 
}
\nonumber  \\
& & \nonumber \\
& = & \left\{\begin{array}{ll}
\displaystyle{(-1)^{a_k} 2^{2\ipart{a_k}} C(v(\beta), \ipart{n} - v(\beta) - \ipart{a_k})}
     & \mbox{if $a_k$ is odd,} \\
     0 & \mbox{if $a_k$ is even.}
\end{array}\right.
\label{E:H-2}
\end{eqnarray}
The last step uses Lemma~\ref{L:G-convolve}.
It follows that $H_-(\alpha) = 0$ if $a_k$ is even. 

Assume from now on
that $a_k$ is odd.
For any composition $\beta=(b_1,\ldots, b_\ell)$, let 
$\ipart{\beta} = \ipart{b_1} + \cdots + \ipart{b_\ell}$.
We will show  that for  $1\le i < k$, 
\begin{equation} \label{E:H-induct}
H_-(\alpha)  =  (-1)^{|\alpha^{i}|} \,
2^{2\ipart{\alpha^{i}}}   \sum_{\beta \ge \alpha_{i} }  
(-1)^{k(\beta) + v(\beta)} C(v(\beta), \ipartn-\ipart{\alpha^{i}} - v(\beta)).
\end{equation}
This holds for $i = k-1$, as can be seen by substituting the first alternative of
\eqref{E:H-2} into the second sum in  \eqref{E:H-1}.
Suppose by induction that \eqref{E:H-induct} is true for some $i$ such that
$k>i>1.$ Since $k(\beta\gamma) = k(\beta) + k(\gamma)$ and $v(\beta\gamma) = 
v(\beta) + v(\gamma)$, we may rewrite \eqref{E:H-induct} as
\begin{eqnarray}
H_-(\alpha) & = & 
(-1)^{|\alpha^i|} 2^{2\ipart{\alpha^i}}
 \sum_{\delta \ge \alpha_{i-1} \atop \gamma \comp a_i}  
 (-1)^{k(\delta\gamma)+v(\delta\gamma)}
 C(v(\delta \gamma), \ipart{n}-\ipart{\alpha^i} - v(\delta \gamma)) 
 \nonumber
\\
& = & (-1)^{|\alpha^i|} 2^{2\ipart{\alpha^i}}
 \sum_{\delta \ge \alpha_{i-1}}  (-1)^{k(\delta)+v(\delta)}
\sum_{\gamma \comp a_i} (-1)^{k(\gamma)+v(\gamma)} 
\nonumber
\\
& & \;\;\;\;\; C(v(\delta) + v(\gamma), \ipart{n}-\ipart{\alpha^i} - v(\delta) - v(\gamma)) 
\nonumber
\\
& = &   (-1)^{|\alpha^{i}|} \,
2^{2\ipart{\alpha^{i}}}  \sum_{\delta \ge \alpha_{i-1}} (-1)^{k(\delta)+v(\delta)} 
\sum_{b = 0}^{\ipart{a_i}}  (-1)^b  
\nonumber 
\\
& & \;\;\;\;\;  C(v(\delta) + b, \ipartn-\ipart{\alpha^i} - v(\delta) - b)
\sum_{\gamma \comp a_i \atop v(\gamma) = b} (-1)^{k(\gamma)}
\nonumber 
\\
& = &  (-1)^{|\alpha^{i-1}|} \,
2^{2\ipart{\alpha^{i}}}  \sum_{\delta \ge \alpha_{i-1}} (-1)^{k(\delta)+v(\delta)}
\label{E:Hminus-induct-step} 
\\
& & \;\;\;\;\;  \sum_{b = 0}^{\ipart{a_i}} C(v(\delta) + b , n/2-\ipart{\alpha^i} - v(\delta) - b) 
 \binom{\ipart{a_i}}{b}. \nonumber
\end{eqnarray}
The last step uses Lemma~\ref{lem:signs}(a).
Now apply
Lemma~\ref{L:G-convolve}, with $i=v(\delta)$ and 
$j = \ipartn - \ipart{\alpha^{i-1}} - v(\delta)$, to simplify the second sum in 
\eqref{E:Hminus-induct-step}. The resulting formula is
$$
H_-(\alpha)  =  (-1)^{|\alpha^{i-1}|} \,
2^{2\ipart{\alpha^{i-1}}}   \sum_{\delta \ge \alpha_{i-1} }  
(-1)^{k(\delta) + v(\delta)} C(v(\delta), \ipartn-\ipart{\alpha^{i-1}} - v(\delta)).
$$
This completes  the proof of \eqref{E:H-induct}.

When $i=1$, \eqref{E:H-induct} becomes
\begin{eqnarray*}
H_-(\alpha) &=& (-1)^{|\alpha^1|} \, 2^{2\ipart{\alpha^1}} \sum_{\beta\comp a_1}
(-1)^{k(\beta) + v(\beta) } C(v(\beta)+1, \ipartn-\ipart{\alpha^1} - v(\beta)) 
\\
& = & (-1)^{|\alpha^1|} \, 2^{2\ipart{\alpha^1}} \sum_{b=0}^{\ipart{a_1}} (-1)^{b+1}
C(b, \ipartn-\ipart{\alpha^1} - b ) 
\sum_{\beta\comp a_1 \atop v(\beta) = b} (-1)^{k(\beta)} \\
& = &
(-1)^{|\alpha^1|} \, 2^{2\ipart{\alpha^1}} \sum_{b=0}^{\ipart{a_1}} 
\binom{\ipart{a_1}}{b} C(b, \ipartn-\ipart{\alpha^1} - b),
\end{eqnarray*}
by Lemma~\ref{lem:signs}(a). Note that
$ |\alpha^1| = n - |a_1| \equiv n-1 \pmod 2$ and
$\ipart{a_1} + \ipart{\alpha^1}  = (n-k_o(\alpha))/2.$ Thus, by
Lemma~\ref{L:G-convolve}, with
$i=1$ and $j=\ipartn-\ipart{a_1} - \ipart{\alpha^1} = \ipart{k_o(\alpha)}$, the 
preceding formula for $H_-(\alpha)$ simplifies to
$$H_-(\alpha) = (-1)^{n-1} 2^{n-k_o(\alpha)} C(0, \ipart{k_o(\alpha)}).$$
\end{case}
\end{proof}


\begin{lemma} \label{L:H+alpha}
Suppose that $n$ is even and $\alpha=(a_1,\ldots, a_k) \comp n$.
Then
$$
H_+(\alpha)=
\left\{\begin{array}{ll}
2^{n-k_o(\alpha)} C(1, k_o(\alpha)/2 -1) 
   & \mbox{if $a_1$ and $a_k$ are odd,} \\
2^{n} 
   & \mbox{if $\alpha = (n)$,} \\
0 
   & \mbox{if $a_1$ or $a_k$ is even and $\alpha\ne (n)$.} \\
\end{array}\right.
$$
\end{lemma}
\begin{proof}
We will consider the cases $\alpha = (n)$ and $\alpha \ne (n)$ separately.
In this proof, the notation $\gamma i$, $\tilde{\gamma}$,
and $\ipart{\beta}$ will have the same meaning as
in the proof of Lemma~\ref{L:H-alpha}.
\setcounter{case}{0}
\begin{case} Suppose that $\alpha = (n)$. 
Using the fact that $\peakaug(\gamma i) =  1+ u(\gamma)$,
a calculation similar to the first part of the proof of Case 1 in Lemma~\ref{L:H-alpha}
yields
\begin{eqnarray*}
H_+((n)) 
                 & = &  C(0,n/2) - \sum_{i=1}^{n-1} \sum_{b=0}^{\ipart{(n-i)}}
                 (-1)^b C(b+1, n/2-b-1) \sum_{\gamma\comp n-i \atop u(\gamma) = b}
                 (-1)^{k(\gamma)}.
\end{eqnarray*}
Make the substitution $j=\ipart{(n-i)}$ and apply
Lemma~\ref{lem:signs}(b) to get
\begin{eqnarray*}
H_+((n)) & = & C(0,n/2) + \sum_{j=0}^{n/2-1} \sum_{b=0}^{j}
                 \binom{j}{b} C(b+1, n/2-b-1)  \\
                 & = & C(0,n/2) + \sum_{b=0}^{n/2-1} C(b+1,n/2-b-1) 
                 \sum_{j=b}^{n/2-1} \binom{j}{b} \\
                 & = & C(0,n/2) + \sum_{b=0}^{n/2-1} C(b+1, n/2-b-1)
                 \binom{n/2}{b+1}.
\end{eqnarray*}
The last step uses the basic identity 
$\sum_{j=b}^{n/2-1} \binom{j}{b} = \binom{n/2}{b+1}$. 
We have shown that $H_+((n)) = \sum_{b=0}^{n/2} \binom{n/2}{b}C(b,n/2-b),$
which equals $2^{n}$ by Lemma~\ref{L:G-convolve}.
\end{case}

\begin{case}
Suppose that $\alpha \ne (n)$.
Using the fact that $\peakaug(\beta\gamma) = 1+ u(\beta) + u(\tilde{\gamma})$
for non-empty compositions $\beta$ and $\gamma$, by calculations similar
to \eqref{E:H-1} and \eqref{E:H-2}, we have
\begin{equation} \label{E:H+base}
H_+(\alpha) = 
(-1)^{a_k} 2^{2\ipart{a_k}} \sum_{\beta\ge \alpha_{k-1}}
(-1)^{k(\beta) + u(\beta) + 1} C(u(\beta)+1, n/2-\ipart{a_k} - u(\beta) -1)
\end{equation}
if $a_k$ is odd, 
and $H_+(\alpha) = 0$ if $a_k$ is even.
%
%
%
%
%

Assume from now on that $a_k$ is odd. 
Note that $u(\delta\gamma) = u(\delta) + v(\gamma)$ for non-empty compositions,
$\delta$ and  $\gamma$. Using \eqref{E:H+base} as the base case,
an inductive calculation similar to the proof of \eqref{E:H-induct} shows that
for $1\le i< k$,
\begin{equation} \label{E:H+induct}
H_+(\alpha)  =  (-1)^{|\alpha^{i}|} \,
2^{2\ipart{\alpha^{i}}}   \sum_{\beta \ge \alpha_{i} }  
(-1)^{k(\beta) + u(\beta) + 1} C(u(\beta)+1, n/2-\ipart{\alpha^{i}} - u(\beta) -1).
\end{equation}
When $i=1$, \eqref{E:H+induct} becomes
\begin{eqnarray*}
H_+(\alpha) &=& (-1)^{|\alpha^1|} \, 2^{2\ipart{\alpha^1}} \sum_{\beta\comp a_1}
(-1)^{k(\beta) + u(\beta) + 1} C(u(\beta)+1, n/2-\ipart{\alpha^1} - u(\beta) -1) 
\\
& = & (-1)^{|\alpha^1|} \, 2^{2\ipart{\alpha^1}} \sum_{b=0}^{\ipart{a_1}} (-1)^{b+1}
C(b+1, n/2-\ipart{\alpha^1} - b -1) 
\sum_{\beta\comp a_1 \atop u(\beta) = b} (-1)^{k(\beta)} \\
& = & \left\{\begin{array}{ll}
\displaystyle{(-1)^{|\alpha^1|} \, 2^{\ipart{\alpha^1}} \sum_{b=0}^{\ipart{a_1}} 
\binom{\ipart{a_1}}{b} C(b+1, n/2-\ipart{\alpha^1} - b -1)}
& \mbox{if $a_1$ is odd,} \\
0 
& \mbox{if $a_1$ is even,}
\end{array}\right.
\end{eqnarray*}
by Lemma~\ref{lem:signs}(b). Note that
$|a_1| + |\alpha^1| = n \equiv 0 \pmod{2}$ and
$\ipart{a_1} + \ipart{\alpha^1}  = (n-k_o(\alpha))/2.$ Thus, by
Lemma~\ref{L:G-convolve}, with
$i=1$ and $j=n/2-\ipart{a_1} - \ipart{\alpha^1} - 1 = k_o(\alpha)/2 - 1$, 
the preceding formula simplifies to
$$H_+(\alpha) = 2^{n-k_o(\alpha)} C(1, k_o(\alpha)/2 -1)$$
when $a_1$ is odd.
\end{case}
The proof is complete.
\end{proof}


\begin{proof}[Proof of Theorem~\ref{T:zetaF}.]
Let $\rho,\psi:\QSym\to\field$ be the linear maps 
defined by the proposed formula for  $\zeta_{+}$ and $\zeta_{-}$, respectively. 
In view of~\eqref{E:FM}, to conclude $\rho = \zeta_+$ and $\psi = \zeta_-$,
it suffices to show that
\begin{equation}  \label{E:zeta-Fproof}
\sum_{\beta \ge \alpha} (-1)^{k(\beta) - k(\alpha)} \psi(F_\beta) = \zeta_-(M_\alpha)
\;\;\;\;\;\mbox{for $\alpha \comp n$;}
\end{equation}
 \begin{equation} \label{E:zeta+Fproof}
\sum_{\beta \ge \alpha} (-1)^{k(\beta) - k(\alpha)} \rho(F_\beta) =\zeta_+(M_\alpha) 
\;\;\;\;\;\mbox{for $\alpha \comp n$, $n$ even.}
\end{equation}

According to \eqref{E:zeta-M}, we can rewrite \eqref{E:zeta-Fproof} as
\begin{eqnarray} 
\lefteqn{
\frac{(-1)^{k(\alpha)}}{2^{2\ipartn}}
\sum_{\beta \ge \alpha} (-1)^{k(\beta) +\peakint(\beta)} 
C(\peakint(\beta), \ipartn - \peakint(\beta)) }  \label{E:zeta-Fproof2}
\\
&  =&  \left\{\begin{array}{ll}
\displaystyle{\frac{(-1)^{k_e(\alpha)}}{2^{2\ipart{k_o(\alpha)}}} C(0,\ipart{k_o(\alpha)})}
& \mbox{if $a_k$ is odd,}
\\
\rule{0pt}{20pt} 0 & \mbox{if $a_k$ is even.} 
\end{array}\right. \nonumber
\end{eqnarray}
Note that the left-hand side of \eqref{E:zeta-Fproof2} equals
$(-1)^{k(\alpha)+1} H_-(\alpha) /2^{2\ipartn}.$
Furthermore, using the facts $k(\alpha) + k_e(\alpha) \equiv k_o(\alpha) \equiv n \pmod 2$
and $\ipart{n} - \ipart{k_o(\alpha)} = (n-k_o(\alpha))/2$, the first alternative
in \eqref{E:zeta-Fproof2} reduces to the equation
$H_-(\alpha) = (-1)^{n-1} 2^{n- k_o(\alpha) } C(0,\ipart{k_o(\alpha)})$; the
second alternative is equivalent to $H_-(\alpha) = 0$.
These formulas for $H_-(\alpha)$ agree with the ones from Lemma~\ref{L:H-alpha}.
This proves \eqref{E:zeta-Fproof2} and hence \eqref{E:zeta-Fproof}.

By similar considerations, \eqref{E:zeta+Fproof} reduces to the formulas
for $H_+(\alpha)$ given in Lemma~\ref{L:H+alpha}.
\end{proof}

\begin{remark}\label{R:odd} The fact that $\zeta_-$ is an odd character, which we know from Theorem~\ref{T:zetaM}, is made transparent by  formula~\eqref{E:zeta-F}. Recall that for any composition $\alpha\comp n$, the antipode of $\QSym$ is given by~\cite[Corollaire 4.20]{Malv}
\begin{equation}\label{E:antipodeF}
S(F_\alpha)=(-1)^nF_{\alpha'}\,,
\end{equation}
where the ribbon diagram of the {\em conjugate} composition $\alpha'$ is obtained by reflecting the
ribbon diagram of $\alpha$ across the line $y=x$. For instance, if $\alpha=(2,3,1,2,2)$ then $\alpha'=(1,2,3,1,2,1)$, as illustrated below:
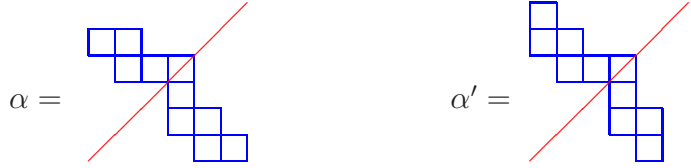
\[ \begin{picture}(60,70)(40,0){\put(-30,20){$\alpha=$}\color{blue}
   \put(0,50){\line(1,0){20}} \put(0,40){\line(1,0){40}}
    \put(10,30){\line(1,0){30}} \put(30,20){\line(1,0){20}}
     \put(30,10){\line(1,0){30}} \put(40,0){\line(1,0){20}}
    \put(0,40){\line(0,1){10}} \put(10,30){\line(0,1){20}}
    \put(20,30){\line(0,1){20}} \put(30,10){\line(0,1){30}}
    \put(40,0){\line(0,1){40}} \put(50,0){\line(0,1){20}}
    \put(60,0){\line(0,1){10}}
    }{\color{red}
    \put(0,0){\line(1,1){60}}
    }\end{picture}
    \qquad\qquad
    \begin{picture}(60,70)(-20,0){\put(-30,20){$\alpha'=$}
    \color{blue}
   \put(50,0){\line(0,1){20}} \put(40,0){\line(0,1){40}}
    \put(30,10){\line(0,1){30}} \put(20,30){\line(0,1){20}}
     \put(10,30){\line(0,1){30}} \put(0,40){\line(0,1){20}}
    \put(40,0){\line(1,0){10}} \put(30,10){\line(1,0){20}}
    \put(30,20){\line(1,0){20}} \put(10,30){\line(1,0){30}}
    \put(0,40){\line(1,0){40}} \put(0,50){\line(1,0){20}}
    \put(0,60){\line(1,0){10}}
    }{\color{red}
    \put(0,0){\line(1,1){60}}
    }\end{picture}\]
It is clear that the number of upper corners in the ribbon diagrams of $\alpha$ and $\alpha'$ are the same, so $\peakint(\alpha)=\peakint(\alpha')$, and hence by~\eqref{E:zeta-F},
\[(\zeta_-)^{-1}(F_\alpha)=(\zeta_-\circ S)(F_\alpha)=(-1)^n\zeta_-(F_{\alpha'})=
(-1)^n\zeta_-(F_\alpha)=\Bar{\zeta}_-(F_\alpha)\,,\]
which shows that $\zeta_-$ is odd. On the other hand, expressing this condition
in the form $\Bar{\zeta}_-\zeta_-=\epsilon$ leads to interesting identities involving bivariate Catalan numbers; see Section~\ref{S:app2}.
\end{remark}

\begin{remark}\label{R:even} Formula~\eqref{E:zeta+F} reveals an analogous property of the character $\zeta_+$, which is not obvious from its definition. 
Recall that if $\alpha=(a_1,\ldots,a_k)$ then $\tilde{\alpha}=(a_k,\ldots,a_1)$ denotes its reversal. The ribbon diagram of $\tilde{\alpha}$ is obtained by reflecting the ribbon diagram of $\alpha$ across the line $y=-x$. Consider the linear map $T:\QSym\to\QSym$ defined by
\begin{equation}\label{E:TF}
T(F_\alpha)=F_{\tilde{\alpha}}\,.
\end{equation}
Since $\alpha\mapsto\tilde{\alpha}$ preserves refinements, the map $T$ is also given by 
\begin{equation}\label{E:TM}
T(M_\alpha)=M_{\tilde{\alpha}}\,.
\end{equation}
{} From either formula it follows easily that $T$ is an antimorphism of coalgebras, a morphism of algebras, and an involution (properties which are shared by the antipode $S$ of $\QSym$).
Moreover, we have
\[\zeta_+\circ T=\zeta_+\,.\]
This property follows from~\eqref{E:zeta+F} and the fact that $\peakaug(\alpha)=\peakaug(\tilde{\alpha})$, which is obvious from~\eqref{E:p+}. Another proof of this fact is given in Proposition~\ref{P:even-odd-T}. 
\end{remark}

\begin{remark}\label{R:even-odd} Let $\alpha=(a_1,\ldots,a_k)$. As mentioned in the preceding remarks,
\[\peakint(\alpha')=\peakint(\alpha) \text{ \ and \ } \peakaug(\tilde{\alpha})=\peakaug(\alpha)\,.\]
We also have the following formulas:
\begin{align}
\label{E:p-rev}
\peakint(\tilde{\alpha}) & =\begin{cases}
\peakint(\alpha) & \text{if $a_1=a_k=1$ or $a_1,a_k \ne 1$,} \\
\peakint(\alpha)-1 & \text{if $a_1\ne 1$ and $a_k=1$,}\\
\peakint(\alpha)+1 & \text{if $a_1= 1$ and $a_k\neq 1$;}
\end{cases}\\
\label{E:p+con}
\peakaug(\alpha') & =\begin{cases}
\peakaug(\alpha) & \text{if exactly one of $a_1$ and $a_k$ is $1$,} \\
\peakaug(\alpha) - 1 & \text{if both $a_1$ and $a_k$ are $1$,} \\
\peakaug(\alpha)+1 & \text{if neither $a_1$ nor $a_k$ is $1$.}
\end{cases}
\end{align}
These formulas follow easily from~\eqref{E:p-} and~\eqref{E:p+}, but they can also be visualized in terms of ribbon diagrams.  We illustrate the second alternative of~\eqref{E:p+con} below, for $\alpha=(1,3,1,2,1)$. Recall that $\peakaug(\alpha)$ is the number of upper corners in the augmented ribbon diagram of $\alpha$.

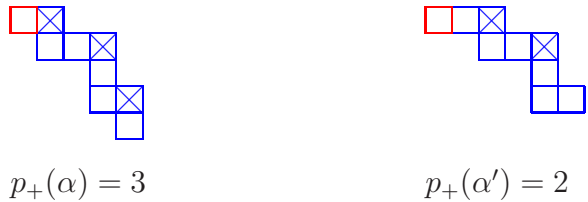
\[ \begin{picture}(60,70)(40,-25){\put(0,-20){$\peakaug(\alpha)=3$}\color{blue}
   \put(10,50){\line(1,0){10}} \put(10,40){\line(1,0){30}}
    \put(10,30){\line(1,0){30}} \put(30,20){\line(1,0){20}}
     \put(30,10){\line(1,0){20}} \put(40,0){\line(1,0){10}}
   \put(10,30){\line(0,1){10}}
    \put(20,30){\line(0,1){20}} \put(30,10){\line(0,1){30}}
    \put(40,0){\line(0,1){40}} \put(50,0){\line(0,1){20}}
     \put(10,50){\line(1,-1){10}} \put(10,40){\line(1,1){10}}
      \put(30,40){\line(1,-1){10}} \put(30,30){\line(1,1){10}}
       \put(40,20){\line(1,-1){10}} \put(40,10){\line(1,1){10}}
    }{\color{red}
    \put(0,50){\line(1,0){10}} \put(0,40){\line(1,0){10}}
    \put(0,40){\line(0,1){10}} \put(10,40){\line(0,1){10}}
    }\end{picture}
    \qquad\qquad
    \begin{picture}(60,70)(-20,-25){\put(-10,-20){$\peakaug(\alpha')=2$}
    \color{blue}
   \put(50,10){\line(0,1){10}} \put(40,10){\line(0,1){30}}
    \put(30,10){\line(0,1){30}} \put(20,30){\line(0,1){20}}
     \put(10,30){\line(0,1){20}} 
    \put(30,10){\line(1,0){20}}
    \put(30,20){\line(1,0){20}} \put(10,30){\line(1,0){30}}
    \put(0,40){\line(1,0){40}} \put(0,50){\line(1,0){20}}
    \put(40,30){\line(-1,1){10}} \put(30,30){\line(1,1){10}}
       \put(20,40){\line(-1,1){10}} \put(10,40){\line(1,1){10}}
    }{\color{red}
    \put(-10,40){\line(1,0){10}}\put(-10,50){\line(1,0){10}}
     \put(-10,40){\line(0,1){10}}\put(0,40){\line(0,1){10}}
    }\end{picture}\]

The statistics $\peakint(\tilde{\alpha})$ and $\peakaug(\alpha')$ enter in some of the formulas in Section~\ref{S:inverses}.
\end{remark}

\section{Application: Identities for bivariate Catalan numbers}\label{S:app2}

The first applications we propose stem from evaluating products of characters
on the basis $F_\alpha$ of $\QSym$. This requires knowledge of the coproduct of $\QSym$ on this basis. Consider all ways of cutting the ribbon diagram of $\alpha$ into two pieces along the common boundary of two squares.  
We include the two trivial cuts along the first and last edges. Thus, if $\alpha$ is a composition of $n$,  there are $n+1$ ways of cutting its
diagram.  For $\alpha=(2,3,1,2,2)$, two non-trivial cuts are shown below.
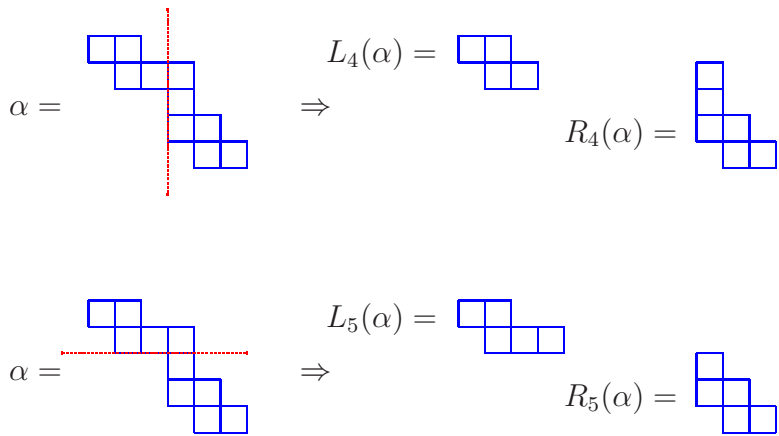
\[ \begin{picture}(60,80)(40,-10){\put(-30,20){$\alpha=$}\color{blue}
   \put(0,50){\line(1,0){20}} \put(0,40){\line(1,0){40}}
    \put(10,30){\line(1,0){30}} \put(30,20){\line(1,0){20}}
     \put(30,10){\line(1,0){30}} \put(40,0){\line(1,0){20}}
    \put(0,40){\line(0,1){10}} \put(10,30){\line(0,1){20}}
    \put(20,30){\line(0,1){20}} \put(30,10){\line(0,1){30}}
    \put(40,0){\line(0,1){40}} \put(50,0){\line(0,1){20}}
    \put(60,0){\line(0,1){10}}
    }{\color{red}
    \put(30,-10){\dashbox{.8}(0,70){}}
    }\put(80,20){$\Rightarrow$}
    \end{picture}
 \begin{picture}(60,60)(-40,-10){\put(-50,40){$L_4(\alpha)=$}\color{blue}
   \put(0,50){\line(1,0){20}} \put(0,40){\line(1,0){30}}
    \put(10,30){\line(1,0){20}} 
     \put(0,40){\line(0,1){10}} \put(10,30){\line(0,1){20}}
     \put(20,30){\line(0,1){20}} \put(30,30){\line(0,1){10}}
     } \end{picture}
 \begin{picture}(60,60)(-40,-10){\put(-20,10){$R_4(\alpha)=$}\color{blue}
    \put(30,20){\line(1,0){20}}
     \put(30,10){\line(1,0){30}} \put(40,0){\line(1,0){20}}
     \put(30,40){\line(1,0){10}}  \put(30,30){\line(1,0){10}}
    \put(40,0){\line(0,1){40}} \put(50,0){\line(0,1){20}}
    \put(60,0){\line(0,1){10}}\put(30,10){\line(0,1){30}}
    }\end{picture}
    \]
    
    \[ \begin{picture}(60,80)(40,-10){\put(-30,20){$\alpha=$}\color{blue}
   \put(0,50){\line(1,0){20}} \put(0,40){\line(1,0){40}}
    \put(10,30){\line(1,0){30}} \put(30,20){\line(1,0){20}}
     \put(30,10){\line(1,0){30}} \put(40,0){\line(1,0){20}}
    \put(0,40){\line(0,1){10}} \put(10,30){\line(0,1){20}}
    \put(20,30){\line(0,1){20}} \put(30,10){\line(0,1){30}}
    \put(40,0){\line(0,1){40}} \put(50,0){\line(0,1){20}}
    \put(60,0){\line(0,1){10}}
    }{\color{red}
    \put(-10,30){\dashbox{.8}(70,0){}}
    }\put(80,20){$\Rightarrow$}
    \end{picture}
 \begin{picture}(60,60)(-40,-10){\put(-50,40){$L_5(\alpha)=$}\color{blue}
   \put(0,50){\line(1,0){20}} \put(0,40){\line(1,0){40}}
    \put(10,30){\line(1,0){30}} 
     \put(0,40){\line(0,1){10}} \put(10,30){\line(0,1){20}}
     \put(20,30){\line(0,1){20}} \put(30,30){\line(0,1){10}}
      \put(40,30){\line(0,1){10}}
     } \end{picture}
 \begin{picture}(60,60)(-40,-10){\put(-20,10){$R_5(\alpha)=$}\color{blue}
    \put(30,20){\line(1,0){20}}
     \put(30,10){\line(1,0){30}} \put(40,0){\line(1,0){20}} 
     \put(30,30){\line(1,0){10}}
    \put(40,0){\line(0,1){30}} \put(50,0){\line(0,1){20}}
    \put(60,0){\line(0,1){10}}\put(30,10){\line(0,1){20}}
    }\end{picture}
    \]
Label the edges between adjacent squares from $0$ to $n$ from left to right
and top to bottom, including the first and last edges. 
Let $L_i(\alpha)$ and $R_i(\alpha)$ be the compositions whose ribbon diagrams are the resulting
pieces after cutting the ribbon diagram of $\alpha$ along edge $i$, with $L_i(\alpha)$ corresponding to the left piece and $R_i(\alpha)$ to the right one.
The two trivial cuts result in $L_0(\alpha)=(\,)$, $R_0(\alpha)=\alpha$, and $L_n(\alpha)=\alpha$, $R_n(\alpha)=(\,)$.  The coproduct of $\QSym$ is~\cite[Corollaire 4.17]{Malv}
\begin{equation}\label{E:coproduct-F}
\Delta(F_\alpha)=\sum_{i=0}^n F_{L_i(\alpha)}\otimes F_{R_i(\alpha)}\,.
\end{equation}
The counit is
\[\epsilon(F_\alpha)=\begin{cases} 1 & \text{ if $\alpha=(\,)$,}\\
0 & \text{ otherwise.} \end{cases}\]

We abbreviate
\begin{align*}
\ell_+^i(\alpha):=p_+\bigl(L_i(\alpha)\bigr)\,, & & r_+^i(\alpha):=p_+\bigl(R_i(\alpha)\bigr)\,,\\
\ell_-^i(\alpha):=p_-\bigl(L_i(\alpha)\bigr)\,, & & r_-^i(\alpha):=p_-\bigl(R_i(\alpha)\bigr)\,.\\
\end{align*}
\begin{proposition}\label{P:app-F}
For any composition $\alpha$ of a positive integer $n$,
\begin{multline}
\label{E:app-F1}
\sum_{j=0}^{\ipartn} (-1)^{\ell^{2j}_+(\alpha)+r_-^{2j}(\alpha)}
C\bigl(\ell_+^{2j}(\alpha),j-\ell_+^{2j}(\alpha)\bigr)
C\bigl(r_-^{2j}(\alpha),\ipartn-j-r_-^{2j}(\alpha)\bigr)\\
=\begin{cases} 4^{\ipart{n}} & \text{ if }\alpha=(n),\\
0 & \text{ otherwise.}
\end{cases}
\end{multline}
\begin{equation}
\label{E:app-F2}
\sum_{i=0}^{n} \frac{(-1)^{\ell^{i}_-(\alpha)+r_-^{i}(\alpha)+i}}
{4^{\ipartfrac{i}+\ipartfrac{n-i}}}
C\bigl(\ell_-^{i}(\alpha),\ipart{i}-\ell_-^{i}(\alpha)\bigr)
C\bigl(r_-^{i}(\alpha),\ipart{(n-i)}-r_-^{i}(\alpha)\bigr)=0\,.
\end{equation}
\end{proposition}
\begin{proof} Since $\zeta_+$ and $\zeta_-$ are the even and odd parts of $\zeta$, we have
\[\zeta(F_\alpha)=(\zeta_+\zeta_-)(F_\alpha)=\sum_{i=0}^n \zeta_+\bigl(F_{L_i(\alpha)}\bigr)\zeta_-\bigl(F_{R_i(\alpha)}\bigr)\,.\]
Note that $\abs{L_i(\alpha)}=i$. Hence, by~\eqref{E:zeta+F}, only
the terms corresponding to even $i$ contribute to this sum. We evaluate these terms using~\eqref{E:zeta-F}, \eqref{E:zeta+F}, and~\eqref{E:zeta-QSym}.
The power of $2$ that results in the denominator is constant and equal to $2\ipartn$. Identity~\eqref{E:app-F1} follows.

Identity~\eqref{E:app-F2} follows similarly from $(\Bar{\zeta}_-\zeta_-)(F_\alpha)=\epsilon(F_\alpha)=0$.
\end{proof}

Various identities for bivariate Catalan numbers may be obtained through special choices of $\alpha$ in Proposition~\ref{P:app-F}. We discuss one
 case involving {\em central Catalan numbers} that is particularly nice.

We refer to the following as central Catalan numbers:
\begin{align}
\label{E:CGnumbers}
C_1(h):=\frac{1}{2}C(2h+1,h+1)\,, & &C_2(h):=\frac{1}{2}C(2h,h+1)\,,\\
C_3(h):=\frac{1}{2}C(2h,h)\,, &  &C_4(h):=\frac{1}{2}C(2h+1,h)\,.\notag
\end{align}
It follows from~\eqref{E:Gnumbers} that
\begin{align*}
C_1(h)=\binom{4h+2}{h}\,, & & C_2(h)=\frac{2h+1}{4h+1}\binom{4h+1}{h}\,,\\
C_3(h)=\frac{1}{2}\binom{4h}{h}\,, & & C_4(h)=\binom{4h+1}{h}\,.
\end{align*}
We have the following wonderful convolution formulas for central Catalan numbers. Note that on each side of each equation, the total sum of the subindices $r$ in $C_r$ is constant (equal to $6$, $7$, and $8$ in each case).
Ira Gessel has shown us how one may derive these identities via generating functions~\cite{Ges-pri}.

\begin{corollary}\label{C:CGnumbers} For any positive integer $h$,
\begin{align}
\label{E:CG6}
\sum_{j=0}^h C_3(j)C_3(h-j) & =2\sum_{j=0}^{h-1}C_2(j)C_1(h-1-j)\,,\\
\label{E:CG7}
\sum_{j=0}^h C_3(j)C_4(h-j) & =\sum_{j=0}^{h}C_2(j)C_3(h-j)+\sum_{j=0}^{h-1}C_1(j)C_1(h-1-j)\,,\\
\label{E:CG8}
\sum_{j=0}^h C_4(j)C_4(h-j) & =2\sum_{j=0}^{h}C_3(j)C_1(h-j)\,.
\end{align}
\end{corollary}
\begin{proof} Let $k$ be a positive integer and $\alpha:=(2,1)^k=(2,1,2,1,\ldots,2,1)$. Formula~\eqref{E:coproduct-F} gives
\[\Delta(F_\alpha)=\sum_{i=0}^k F_{(2,1)^i}\otimes F_{(2,1)^{k-i}}+
\sum_{i=0}^{k-1} F_{(2,1)^i,2}\otimes F_{1,(2,1)^{k-1-i}}+
\sum_{i=0}^{k-1} F_{(2,1)^i,1}\otimes F_{1,1,(2,1)^{k-1-i}}\,.\]
Equation~\eqref{E:app-F1} then leads to
\[\sumsub{i=0\\ i \text{ even}}^k C(i,i/2)C(k-i,\ipart{(k-i)})
=\sum_{i=0}^{k-1} C(i,\ipart{i}+1)C(k-1-i,\ipart{k}-\ipart{i})\,.\]
When $k=2h$ this is~\eqref{E:CG6}; when $k=2h+1$ this is~\eqref{E:CG7}.

Similarly, equation~\eqref{E:app-F2} for $\alpha=(2,1)^k$ leads to ~\eqref{E:CG8} when $k=2h$ (and to a trivial identity when $k=2h+1$).
\end{proof}

More identities may be derived from~\eqref{E:zeta-F} and~\eqref{E:zeta+F} by imposing the fact that $\zeta_-$ and $\zeta_+$ are morphisms of algebras.
The multiplication of two basis elements $F_\alpha$ and $F_\beta$ is most easily described in terms of permutations. It is thus convenient to
work on a larger Hopf algebra $\SSym$, of which $\QSym$ is a quotient.
This is the object of the next section. 

\section{Identities for bivariate Catalan numbers via the Hopf algebra of permutations}
\label{S:permutations}

The Hopf algebra of permutations $\SSym$ has a linear basis $\{F_\sigma\}$ indexed by permutations $\sigma\in S_n$, $n\geq 0$. The multiplication of two basis elements is as follows. Given $\sigma\in S_n$ and $\tau\in S_m$, let
$\ISh(\sigma,\tau)$ be the set of all shuffles of the words $\sigma(1),\ldots,\sigma(n)$ and $n+\tau(1),\ldots,n+\tau(m)$. Then
\begin{equation}\label{E:prod-F}
F_\sigma\cdot F_\tau=\sum_{\rho\in\ISh(\sigma,\tau)}F_\rho\,.
\end{equation}
For example,
 \begin{eqnarray*}
  F_{\blue{12}}\cdot F_{\red{312}} &=& \ \ \ 
    F_{\blue{12}\red{534}}\,+\,
    F_{\blue{1}\red{5}\blue{2}\red{34}}\,
    +\,F_{\blue{1}\red{53}\blue{2}\red{4}}\,+\,
    F_{\blue{1}\red{534}\blue{2}}\,+\,
    F_{\red{5}\blue{12}\red{34}}\\
    && +\,F_{\red{5}\blue{1}\red{3}\blue{2}\red{4}}\,+\,
    F_{\red{5}\blue{1}\red{34}\blue{2}}\,+\,
    F_{\red{53}\blue{12}\red{4}}\,
    +\,F_{\red{53}\blue{1}\red{4}\blue{2}}\,+\,
    F_{\red{534}\blue{12}}\,.
 \end{eqnarray*}

For more information on the Hopf algebra structure of $\SSym$ see~\cite{AS}.

The {\em descent set} of a permutation $\sigma\in S_n$ is
\[ \Des(\sigma):=\{ i \in [n-1] \mid \sigma(i) > \sigma(i+1)\}\,.\]
Let $D(\sigma)= (a_1, \ldots, a_k)$ be the composition of $n$ such that
$\{a_1, a_1+a_2, \ldots, a_1+\cdots + a_{k-1}\}=\Des(\sigma)$. The map
\begin{equation} \label{E:descentmap}
   \begin{array}{rcrcl}
     \calD &:& \SSym&\longrightarrow& \QSym\\
           & &F_\sigma&\longmapsto& F_{D(\sigma)}\rule{0pt}{14pt}
   \end{array}
 \end{equation}
is a surjective morphism of graded Hopf algebras~\cite[Th\'eor\`emes 5.12, 5.13, and 5.18]{Malv}.

Let $\zetaS$ denote the pull-back of the universal character $\zeta$ of $\QSym$ via the morphism $\calD$:
\[\zetaS:=\zeta\circ\calD\,.\]
According to Lemma~\ref{L:even-odd-functorial}, the even and odd parts of $\zetaS$ are
\begin{equation}\label{E:zetaSQ}
(\zetaS)_+=\zeta_+\circ\calD \text{ \ and \ } (\zetaS)_-=\zeta_-\circ\calD\,.
\end{equation}
We describe these characters directly in terms of permutations.

As in~\cite{ABN}, consider the following two slightly differing notions of a {\em peak set} of a permutation $\sigma\in S_n$:
\begin{align*}
\Peakint(\sigma) &:= \{ i \in [n-1] \mid i\neq 1,\ \sigma(i-1) < \sigma(i) > \sigma(i+1)\}\,,\\
\Peakaug(\sigma) &:= \{ i \in [n-1] \mid \sigma(i-1) < \sigma(i) > \sigma(i+1)\}\,,
\end{align*}
where we agree that $\sigma(0) = 0$. For instance, if $\sigma=312546$ then
$\Peakint(\sigma) = \{4\} $ and $\Peakaug(\sigma) = \{1,4\}$. The study of peak enumeration has a long history, but the connections between peaks and quasi-symmetric functions originate in work of Stembridge~\cite{Ste97}.

\[ \begin{picture}(120,90)(0,-30){\color{blue}
   \put(40,0){\line(1,1){40}}\put(100,20){\line(1,1){20}} 
    \put(20,20){\line(1,-1){20}}\put(80,40){\line(1,-1){20}}
   \put(20,20){\circle*{3}}\put(40,0){\circle*{3}}
    \put(60,20){\circle*{3}}\put(80,40){\circle*{3}}\put(100,20){\circle*{3}}
    \put(120,40){\circle*{3}}
     \put(17,25){$3$}\put(37,5){$1$} \put(55,25){$2$}\put(77,45){$5$}
     \put(97,25){$4$}\put(117,45){$6$}
   }{\color{red} \put(0,0){\circle*{3}} \put(0,0){\line(1,1){20}}
    \put(-3,5){$0$}
     }\put(-3,-20){$0$}\put(17,-20){$1$}\put(37,-20){$2$}\put(55,-20){$3$}
     \put(77,-20){$4$}\put(97,-20){$5$}\put(117,-20){$6$}
     \put(0,-7){\circle*{3}} \put(20,-7){\circle*{3}}\put(40,-7){\circle*{3}}
    \put(60,-7){\circle*{3}}\put(80,-7){\circle*{3}}\put(100,-7){\circle*{3}}
    \put(120,-7){\circle*{3}}   \put(0,-7){\line(1,0){120}}
       \end{picture}\]

Note that $\Peakint(\sigma)$ and $\Peakaug(\sigma)$ depend only on  $\Des(\sigma)$. In fact, 
\begin{align*}
i\in \Peakaug(\sigma) & \iff i\in\Des(\sigma) \text{ and }
i-1\notin \Des(\sigma)\,,\\
i\in\Peakint(\sigma) & \iff i\in\Peakaug(\sigma) \text{ and }i\neq 1\,.
\end{align*}
We write $\peakint(\sigma) := \#\Peakint(\sigma)$ and $\peakaug(\sigma) := \#\Peakaug(\sigma)$. Let $\alpha:=D(\sigma)$. It follows from~\eqref{E:p-}, ~\eqref{E:p+}, and the above that
\begin{equation}\label{E:peaknumbers}
\peakint(\sigma) = \peakint(\alpha)\;\;\;\mbox{ and } \;\;\;\peakaug(\sigma)=\peakaug(\alpha)\,.
\end{equation}

\begin{proposition}\label{P:zetaS} Let $\sigma\in S_n$ be a permutation. Then
\begin{eqnarray}
(\zetaS)_{-}(F_\sigma) & = & 
\displaystyle{\frac{(-1)^{\peakint(\sigma)}}{2^{2\ipartn}} 
C\left(\peakint(\sigma),\lfloor n/2 \rfloor-\peakint(\sigma)\right).} 
\label{E:zetaS-}\\
(\zetaS)_{+}(F_\sigma) & = &  
\left\{\begin{array}{ll}
\displaystyle{\frac{(-1)^{\peakaug(\sigma)}}{2^{n}} 
C\left(\peakaug(\sigma),n/2 -\peakaug(\sigma)\right)}
& \mbox{ if $n$ is even,} \\
0 & \mbox{ if $n$ is odd.} \end{array} \right.
\label{E:zetaS+}
\end{eqnarray}
\end{proposition}
\begin{proof} This follows at once from Theorem~\ref{T:zetaF} together with~\eqref{E:descentmap},~\eqref{E:zetaSQ} and~\eqref{E:peaknumbers}.
\end{proof}

We may now deduce the identities for  bivariate Catalan numbers announced at the end of Section~\ref{S:app2}.

\begin{corollary}\label{C:allperms} For any non-negative integer $n$,
\begin{align}
\label{E:allperms-}
 \sum_{\sigma \in S_n} (-1)^{\peakint(\sigma)} C\bigl(\peakint(\sigma), \ipartn - \peakint(\sigma)\bigr) & =2^{2\ipartn }\,,\\
\intertext{ and for any even positive integer $n$,}
\label{E:allperms+}
 \sum_{\sigma \in S_n} (-1)^{\peakaug(\sigma)} C\bigl(\peakaug(\sigma), n/2 - \peakaug(\sigma)\bigr) & =0\,.
\end{align}
\end{corollary}
\begin{proof} Consider the $n$-th power of the basis element indexed by the permutation $1\in S_1$ in the algebra $\SSym$.  Since $\ISh(1,1,\ldots,1)=S_n$, ~\eqref{E:prod-F} gives
\[(F_1)^n=\sum_{\sigma\in S_n}F_\sigma\,.\]
Applying $(\zetaS)_-$ to both sides, using that it is a morphism of algebras and~\eqref{E:zetaS-}, we deduce~\eqref{E:allperms-}. Similarly, applying $(\zetaS)_+$ we obtain~\eqref{E:allperms+}.
\end{proof}

Let $1_n$ denote the identity permutation in $S_n$. Note that $\ISh(1_n,1_m)$ is the set of  those permutations $\sigma\in S_{n+m}$ such that
\[\sigma^{-1}(1)<\cdots<\sigma^{-1}(n),\ \sigma^{-1}(n+1)<\cdots<\sigma^{-1}(n+m)\,.\]

\begin{corollary}\label{C:shuffles} Let $n,m$ be non-negative integers. Then
\begin{multline}
\label{E:shuffle-}
\sum_{\sigma \in \ISh(1_n,1_m)} (-1)^{\peakint(\sigma)} C\bigl(\peakint(\sigma), \ipart{(n+m)} -\peakint(\sigma)\bigr) \\
 =\begin{cases} C(0, \ipartn) C(0, \ipart{m}) & \text{ if $n$ or $m$ is even,}\\
 4C(0, \ipartn) C(0, \ipart{m}) & \text{ if $n$ and $m$ are odd.}\end{cases}
 \end{multline}
 \begin{multline}
  \label{E:shuffle+}
\sum_{\sigma \in \ISh(1_n,1_m)} (-1)^{\peakaug(\sigma)} C\bigl(\peakaug(\sigma), (n+m)/2 -\peakaug(\sigma)\bigr) \\
=\begin{cases} C(0, n/2) C(0, m/2)& \text{ if $n$ and $m$ are even,}\\
 0 & \text{ if $n$ and $m$ are odd.}\end{cases}
 \end{multline}
\end{corollary}
\begin{proof} By~\eqref{E:prod-F},
\[F_{1_n}\cdot F_{1_m}=\sum_{\sigma \in \ISh(1_n,1_m)} F_\sigma\,.\]
Applying $(\zetaS)_-$ we deduce~\eqref{E:shuffle-} and applying $(\zetaS)_+$ we deduce~\eqref{E:shuffle+}.
\end{proof}

\section{Inverses of canonical characters and more applications}\label{S:inverses}

The set of characters of a Hopf algebra $\H$ is a group under the convolution product (Section~\ref{S:intro}). The inverse of a character $\varphi$ is $\varphi\circ S$,
where $S$ is the antipode of $\H$.

\begin{proposition}\label{P:zetainv} The inverse of
the universal character $\zeta$ of $\QSym$ is explicitly given by
\begin{equation}\label{E:zetainv-QSym}
\zeta^{-1}(M_\alpha)=(-1)^{k(\alpha)}\,,\qquad
\zeta^{-1}(F_\alpha)=\begin{cases}
(-1)^{\abs{\alpha}} & \text{ if $\alpha=(1,1,\ldots,1)$,} \\
  0       & \text{ if not.}
\end{cases}
\end{equation}
\end{proposition}
\begin{proof} This follows at once from the explicit formulas for the antipode of $\QSym$~\eqref{E:antipodeM} and~\eqref{E:antipodeF}.
\end{proof}

Inverting the canonical decomposition  $\zeta=\zeta_+\zeta_-$ we obtain
\[\zeta^{-1}=(\zeta_-)^{-1}(\zeta_+)^{-1}=(\zeta_+)^{-1}\Bigl(\zeta_+(\zeta_-)^{-1}(\zeta_+)^{-1}\Bigr)=(\zeta_+)^{-1}\Bigl(\zeta_+\Bar{\zeta}_-(\zeta_+)^{-1}\Bigr)\,.\]

The set of even characters is a subgroup of the group of characters, and the set of odd characters is closed under conjugation by even characters~\cite[Proposition 1.7]{ABS}. In particular, $(\zeta_+)^{-1}$ is even and $\zeta_+\Bar{\zeta}_-(\zeta_+)^{-1}$ is odd. According to Lemma~\ref{L:even-odd-decomposition}, these are the even and odd parts of $\zeta^{-1}$:
\begin{equation}\label{E:even-odd-zetainv}
(\zeta^{-1})_+=(\zeta_+)^{-1} \text{ \ and \ } (\zeta^{-1})_-=\zeta_+\Bar{\zeta}_-(\zeta_+)^{-1}\,.
\end{equation}
We provide explicit descriptions for these characters below. First, we analyze the behavior of the map $T:\QSym\to\QSym$ (Remark~\ref{R:even}) with respect to the canonical decomposition of $\zeta$.

\begin{proposition}\label{P:even-odd-T} We have
\begin{equation}\label{E:even-odd-T}
\zeta_+\circ T=\zeta_+ \text{ \ and \ }\zeta_-\circ T=\bigl((\zeta^{-1})_-\bigr)^{-1}\,.
\end{equation}
\end{proposition}
\begin{proof} Both $S$ and $T$ are antimorphism of coalgebras and morphisms of algebras $\QSym\to\QSym$. 
In addition, $\zeta\circ S=\zeta^{-1}$ and $\zeta\circ T=\zeta$; the latter being an immediate consequence of~\eqref{E:zeta-QSym} and~\eqref{E:TF}. Therefore, $T\circ S:\QSym\to\QSym$ is a morphism of Hopf algebras such that
$\zeta\circ T\circ S=\zeta^{-1}$. According to Lemma~\ref{L:even-odd-functorial} we have
\[(\zeta^{-1})_+=\zeta_+\circ T\circ S \text{ \ and \ } (\zeta^{-1})_-=\zeta_-\circ T\circ S\,.\]
Composing with $S$ we find
$\zeta_+\circ T=(\zeta^{-1})_+\circ S=\bigl((\zeta^{-1})_+\bigr)^{-1}=\zeta_+$, by~\eqref{E:even-odd-zetainv}, and $\zeta_-\circ T=(\zeta^{-1})_-\circ S=\bigl((\zeta^{-1})_-\bigr)^{-1}$.
\end{proof}

\begin{remark}\label{R:even-odd-T}
The fact that $\zeta_+\circ T=\zeta_+$ was observed in Remark~\ref{R:even}.
The fact that $\zeta_-\circ T=\bigl((\zeta^{-1})_-\bigr)^{-1}$ may be rewritten as follows:
\begin{equation}\label{E:odd-T}
(\zeta^{-1})_-=\Bar{\zeta}_-\circ T\,.
\end{equation}
Indeed, any morphism (or antimorphism) preserves inverses, so 
\[(\zeta^{-1})_-=\bigl(\zeta_-\circ T\bigr)^{-1}=(\zeta_-)^{-1}\circ T=\Bar{\zeta}_-\circ T\,,\]
since $\zeta_-$ is odd. We use this formula  to calculate $(\zeta^{-1})_-$ below.
\end{remark}

\begin{theorem}\label{T:zetainvM}
Let  $\alpha=(a_1,\ldots,a_k)$ be a composition of a non-negative  integer $n$. Then
\begin{align}
\label{E:zetainv-M}
(\zeta^{-1})_-(M_\alpha)  & = \begin{cases}
{\displaystyle \frac{(-1)^{k(\alpha)} }{2^{2\ipart{k_o(\alpha)}}} C\bigl(0,\ipart{k_o(\alpha)}\bigr)} & \text{ if $a_1$ is odd,} \\
\rule{0pt}{20pt}  0   & \text{ if $a_1$ is even.}
\end{cases}\\\label{E:zetainv+M}
(\zeta^{-1})_+(M_\alpha)  & = \begin{cases}
{\displaystyle \frac{(-1)^{k(\alpha)}}{2^{k_o(\alpha)}} C\bigl(0,k_o(\alpha)/2\bigr)} & \text{ if $n$ is even,} \\
\rule{0pt}{20pt}  0       & \text{ if $n$ is odd;}
\end{cases}
\end{align}
\end{theorem}
\begin{proof} According to~\eqref{E:odd-T} and~\eqref{E:TM},
\[(\zeta^{-1})_-(M_{\alpha})=\Bar{\zeta}_-\bigl( T(M_{\alpha})\bigr)=(-1)^n\zeta_-(M_{\tilde{\alpha}})\,,\]
which we can evaluate with~\eqref{E:zeta-M}. To see that this results in~\eqref{E:zetainv-M}, note that
since $\tilde{\alpha}$ is the reversal of $\alpha$, the last part of $\tilde{\alpha}$
is $a_1$, $k_o(\tilde{\alpha})=k_o(\alpha)$, and  by~\eqref{E:oddparts}
\[n+k_e(\tilde{\alpha})=n+k_e(\alpha)\equiv k_o(\alpha)+k_e(\alpha)=k(\alpha)\,.\] 


To settle the remaining identity we give a direct argument. Let $\varphi :\QSym\to\field$ be the linear functional defined by~\eqref{E:zetainv+M}. We show that $\zeta_+ \varphi = \epsilon$, which implies $\varphi=(\zeta_+)^{-1}=(\zeta^{-1})_+$.

Since  $\varphi(M_{(\,)}) = 1$, we have $(\zeta_+ \varphi)(M_{(\,)}) =  1 = 
\epsilon(M_{(\,)})$. 

Assume from now on that $n > 0$. We need to show  $(\zeta_+ \varphi)(M_\alpha) = 0$. Write $\alpha=(a_1,\ldots, a_k)$ and recall the notations $\alpha_i$ and $\alpha^i$ from~\eqref{E:convolution-QSym}.

If $n$ is odd then for each $i$ one of $\abs{\alpha_i}$ and $\abs{\alpha^i}$ is odd, so
every term in the expansion~\eqref{E:convolution-QSym} of $(\zeta_+ \varphi)(M_\alpha)$ is $0$, by~\eqref{E:zeta+M} and~\eqref{E:zetainv+M}. Assume from now on that $n$ is even.

If $k = 1$ (i.e., $\alpha = (n)$), then  $(\zeta_+ \varphi)(M_{(n)}) =
\varphi(M_{(n)}) + \zeta_+(M_{(n)}) = -1 + 1  = 0.$

If $k > 1$ and $a_1$ is even, 
then $\varphi(M_{\alpha}) = - \varphi(M_{\alpha^1})$ (since $\alpha^1$ has one less part than $\alpha$ and the same number of odd parts),
and $\zeta_+(M_{\alpha_i}) = 0$ when $i > 1$ by~\eqref{E:zeta+M}.
Thus,
$$(\zeta_+ \varphi)(M_\alpha) = 
\varphi(M_\alpha) + \zeta_+(M_{(a_1)}) \varphi(M_{\alpha^1}) = 
\varphi(M_\alpha) + \varphi(M_{\alpha^1}) =0.$$

Suppose finally that $k>1$ and $a_1$ is odd.
Using~\eqref{E:convolution-QSym} and~\eqref{E:zeta+M} we compute
\begin{eqnarray*}
(\zeta_+ \varphi)(M_\alpha) 
& = & 
\frac{(-1)^{k(\alpha)}}{2^{k_o(\alpha)}} B\bigl(k_o(\alpha)/2\bigr) + 
\sum_{{1\le  i \le k \atop a_i  \text{ odd}} \atop \abs{\alpha_i}  \text{ even}}
\frac{(-1)^{k_e(\alpha_i)+1}}{2^{k_o(\alpha_i)-1}} C\bigl(k_o(\alpha_i)/2-1\bigr)
\frac{(-1)^{k(\alpha^i)}}{2^{k_o(\alpha^i)}} B\bigl(k_o(\alpha^i)/2\bigr)
\\
& = & 
\frac{(-1)^k}{2^{k_o(\alpha)}} B\bigl(k_o(\alpha)/2\bigr) +
\frac{(-1)^{k+1}}{2^{k_o(\alpha) - 1}}
\sum_{{1\le  i \le k \atop a_i  \text{ odd}} \atop \abs{\alpha_i}  \text{ even}}
C\bigl(k_o(\alpha_i)/2-1\bigr) B\bigl(k_o(\alpha^i)/2\bigr)\,. 
\end{eqnarray*}
To combine the signs in the last step we used
$k=k(\alpha) \equiv k_e(\alpha_i) + k(\alpha^i) \mod 2$, which holds by~\eqref{E:oddparts} and since $\abs{\alpha_i}$, $\abs{\alpha^i}$, and $\abs{\alpha}$ are all even.
The argument may now be completed as in the proof of Theorem~\ref{T:zetaM}, Case~\ref{case5}. 
\end{proof}

We present some applications. The first one is an identity which appears in~\cite[Section 4.2, Example 2]{Rio} (near the bottom of page 130).

\begin{corollary}\label{C:app-zetainv-M} For any positive integer $m$,
\begin{equation}\label{E:app-zetainv-M}
\sum_{j=0}^{m-1}2^{2m-2j-1}C(j)=2^{2m}-\binom{2m}{m}\,.
\end{equation}
\end{corollary}
\begin{proof} We have $\zeta^{-1}=(\zeta^{-1})_+(\zeta^{-1})_-=(\zeta_+)^{-1}(\zeta^{-1})_-$, so $(\zeta^{-1})_-=\zeta_+\zeta^{-1}$. We evaluate both sides on $M_{\alpha}$ with $\alpha=(1,1,\ldots,1)\comp 2m$. Equations~\eqref{E:zeta+M},~\eqref{E:zetainv-QSym}, and~\eqref{E:zetainv-M} lead to
\[\frac{1}{2^{2m}} C(0,m)=1+\sumsub{i=2 \\ i \text{ even}}^{2m}\frac{-1}{2^{i}} C(1,i/2-1)\,.\]
Letting $j=2i+2$ we obtain~\eqref{E:app-zetainv-M}.
\end{proof}

The following identity is analogous to~\eqref{E:preapp-antipodeM}.
\begin{corollary}\label{C:app-zetainv+M} For any composition $\beta$ of an even integer,
\begin{equation}\label{E:app-zetainv+M} 
\sumsub{\alpha\leq\beta\\ a_1,a_k\text{ odd}}(-1)^{k_e(\alpha)}2^{k_o(\beta)-k_o(\alpha)+1}C\bigl(k_o(\alpha)/2-1\bigr)=2^{k_o(\beta)}-\binom{k_o(\beta)}{k_o(\beta)/2}\,;
\end{equation}
the sum being over those compositions $\alpha=(a_1,\ldots,a_k)$ whose first and last part are odd and which are refined by $\beta$.
\end{corollary}
\begin{proof} By~\eqref{E:even-odd-zetainv},
\[\zeta_+\circ S=(\zeta_+)^{-1}=(\zeta^{-1})_+\,.\]
Therefore, for any composition $\beta$ of $n$ we have, by~\eqref{E:antipodeM},
\[(-1)^{k(\beta)}\sum_{\alpha\leq\tilde{\beta}}\zeta_+(M_\alpha)=(-1)^n(\zeta^{-1})_+(M_\beta)\,.\]
Hence,  by~\eqref{E:zeta+M} and~\eqref{E:zetainv+M}, and since $n$ is even,
\[(-1)^{k(\beta)}+(-1)^{k(\beta)}\sumsub{\alpha\leq\tilde{\beta}\\ a_1,a_h\text{ odd}}
 \frac{(-1)^{k_e(\alpha)+1}}{2^{k_o(\alpha) }} C(1, k_o(\alpha)/2-1)=
   \frac{(-1)^{k(\beta)}}{2^{k_o(\beta)}} C(0,k_o(\beta)/2)\,.\]
The above sum remains unchanged if we replace $\beta$ for $\tilde{\beta}$.
Multiplying by $2^{k_o(\beta)}$ gives~\eqref{E:app-zetainv+M}.
\end{proof}

\smallskip

On the fundamental basis, the even and odd parts of $\zeta^{-1}$ are most easily described in terms of the statistics $\peakaug(\alpha')$ and
$\peakint(\tilde{\alpha})$ (Remark~\ref{R:even-odd}). They may be described directly in terms of $\peakaug(\alpha)$ and $\peakint(\alpha)$ by means of formulas~\eqref{E:p-rev} and~\eqref{E:p+con}.  
 
\begin{theorem}\label{T:zetainvF}
Let $\alpha = (a_1,\ldots, a_k)$ be a composition of a non-negative  integer $n$. Then
\begin{align}
\label{E:zetainv-F}
(\zeta^{-1})_-(F_\alpha) & = 
{\displaystyle \frac{(-1)^{n+\peakint(\tilde{\alpha})}}{2^{2\ipart{n}}}   
C\bigl(\peakint(\tilde{\alpha}),\ipart{n}-\peakint(\tilde{\alpha})\bigr)  } \\
\label{E:zetainv+F}
(\zeta^{-1})_+(F_\alpha) & = 
\left\{ \begin{array}{ll}
\displaystyle{\frac{(-1)^{\peakaug(\alpha')}}{2^{n}} C\bigl(\peakaug(\alpha'), n/2 - \peakaug(\alpha')\bigr)}
& \text{if $n$ is even,} \\
\rule{0pt}{20pt} 0
& \text{if $n$ is odd.} 
\end{array}\right.
\end{align}
\end{theorem}
\begin{proof}
 Assume that $n$ is even, $n \ne 0$. We have, by~\eqref{E:antipodeF} and~\eqref{E:even-odd-zetainv},
\[(\zeta^{-1})_+(F_\alpha)=(\zeta_+)^{-1}(F_\alpha)=(\zeta_+\circ S)(F_\alpha)=\zeta_+(F_{\alpha'})\,.\]
 Formula~\eqref{E:zetainv+F}  now follows from~\eqref{E:zeta+F}.
  Formula~\eqref{E:zetainv-F}  follows similarly from~\eqref{E:zeta-F} and~\eqref{E:odd-T}.
\end{proof}

As an application we derive a recursive formula for the bivariate Catalan numbers.

 \begin{corollary}\label{C:gessel-rec}
For any non-negative integers $a$, $b$, and $c$,
\begin{equation}\label{E:gessel-rec}
C(b,a+c)= 4^c C(b,a)-\sum_{j=1}^c 4^{c-j} C(b+1,a+j-1)\,.
\end{equation}
\end{corollary}
\begin{proof}We evaluate both sides of $(\zeta^{-1})_-=\zeta_+\zeta^{-1}$ on
$F_{\alpha}$ with \[\alpha=(\underbrace{1,\ldots,1}_{2a+1},\underbrace{2,\ldots,2}_{b},\underbrace{1,\ldots,1}_{2c})\,.\]

According to~\eqref{E:coproduct-F},
 \[(\zeta_+\zeta^{-1})(F_\alpha)=\sum_{i=0}^n \zeta_+(F_{L_i(\alpha)}) \zeta^{-1}(F_{R_i(\alpha)})\,.\]
 Several terms in this expansion vanish, in view of~\eqref{E:zeta+F} and~\eqref{E:zetainv-QSym}: all those for which the size of $L_i(\alpha)$ is odd and all those for which  $R_i(\alpha)\neq(1,\ldots,1)$.
We are left with
\begin{align*}
(\zeta_+\zeta^{-1})(F_\alpha) &=
\zeta_+(F_{(\underbrace{1,\ldots,1}_{2a+1},\underbrace{2,\ldots,2}_{b-1},{\displaystyle 1})})   
\zeta^{-1}(F_{(\underbrace{1,\ldots,1}_{2c+1})})\\
&+\sumsub{i=1\\ i \text{ odd}}^{2c-1}
\zeta_+(F_{(\underbrace{1,\ldots,1}_{2a+1},\underbrace{2,\ldots,2}_{b},\underbrace{1,\ldots,1}_{i})})
\zeta^{-1}(F_{(\underbrace{1,\ldots,1}_{2c-i})})\\
&=\frac{(-1)^{b}}{2^{2(a+b)}}C(b,a)(-1)^{2c+1}+
\sumsub{i=1\\ i \text{ odd}}^{2c-1}\frac{(-1)^{b+1}}{2^{2(a+b)+i+1}}C(b+1,a+(i+1)/2-1)(-1)^{2c-i}\\
&=\frac{(-1)^{b+1}}{2^{2(a+b)}}C(b,a)+
\sum_{j=1}^{c}\frac{(-1)^{b}}{2^{2(a+b+j)}}C(b+1,a+j-1)\,.
\end{align*}
We used that $\peakaug(\underbrace{1,\ldots,1}_{2a+1},\underbrace{2,\ldots,2}_{b-1},1)=b$ and for $i\geq 1$,
$\peakaug(\underbrace{1,\ldots,1}_{2a+1},\underbrace{2,\ldots,2}_{b},\underbrace{1,\ldots,1}_{i})=b+1$.

On the other hand, $\peakint(\alpha)=b$, so by~\eqref{E:zetainv-F},
\[(\zeta^{-1})_-(F_{\alpha})=\frac{(-1)^{2(a+b+c)+1+b}}{2^{2(a+b+c)}}C(b,a+c)\,.\]
Equating $(\zeta^{-1})_-(F_{\alpha})$ to $(\zeta_+\zeta^{-1})(F_\alpha)$ gives
~\eqref{E:gessel-rec}.
\end{proof}

A few special cases of Corollary~\ref{C:gessel-rec} are worth stating.
We obtain formulas expressing  a central binomial coefficient or a Catalan number in terms of  bivariate Catalan numbers.

 \begin{corollary}\label{C:binomial-catalan-gessel}
For any non-negative integers $b$ and $c$,
\begin{align}
\label{E:binomial-gessel}
\binom{2b}{b} &=\frac{1}{4^c}C(b,c)+\sum_{j=1}^c \frac{1}{4^j} C(b+1,j-1)\,,\\
\label{E:catalan-gessel}
2C(b) &=\frac{1}{4^c}C(b,c+1)+\sum_{j=1}^c \frac{1}{4^j} C(b+1,j)\,.
\end{align}
\end{corollary}
\begin{proof} These follow by choosing $a=0$ and $a=1$ in~\eqref{E:gessel-rec}.
\end{proof}
 
  \begin{corollary}\label{C:associator} Let $H(a,b,c):=C(a,b+c)-C(b,a+c)$. Then
  \begin{equation}\label{E:associator}
\frac{1}{4^c}H(a,b,c)=\sum_{j=1}^c\frac{1}{4^j}H(b+1,a+1,j-2)\,.
\end{equation}
  \end{corollary}
\begin{proof} This follows from~\eqref{E:gessel-rec} since $C(a,b)=C(b,a)$.
\end{proof}

\section*{Appendix. Bivariate Catalan numbers as binomial coefficients}\label{A:binomial}

Using the easily verified formula
$$
C(m,n) = (-1)^n 4^{m+n} \binom{m-1/2}{m+n}\,,
$$
one may restate all results in this paper in terms of binomial coefficients.
This allows for simplifications in some of the formulas. 
We list here our main results in this notation. Let $\alpha=(a_1,\ldots,a_k)$ be a composition of $n$.

\begin{align}
\zeta_{-}(M_\alpha) &= \left\{\begin{array}{ll}
{\displaystyle (-1)^{k_e(\alpha)+\ipartfrac{k_o(\alpha)}}
 \binom{-1/2}{ \ipart{k_o(\alpha)}} } & 
  \mbox{if $a_k$ is odd,}  \\
0 & \mbox{if $a_k$ is even;}
\end{array} \right. 
\\
\zeta_{+}(M_\alpha)  &=  \left\{\begin{array}{ll}
\rule{0pt}{20pt}{\displaystyle (-1)^{k_e(\alpha)+\frac{k_o(\alpha)}{2}}
 \binom{1/2}{ k_o(\alpha)/2} }& \mbox{if $a_1$ and $a_k$ are odd and $n$ is even,} \\
\rule{0pt}{20pt} 1 & \mbox{if $\alpha=(n)$ and $n$ is even,}\\
0 & \mbox{otherwise;}  \end{array}\right. \\
\rule{0pt}{20pt} \zeta_{-}(F_\alpha) & =  
(-1)^{\ipartfrac{n}}\binom{\peakint(\alpha)-1/2}{\ipartn}\,; \\
\zeta_{+}(F_\alpha) & =   
 \left\{\begin{array}{ll}
\rule{0pt}{20pt}{\displaystyle (-1)^{\frac{n}{2}}\binom{\peakaug(\alpha)-1/2}{n/2} }& \mbox{ if $n$ is even,} \\
\rule{0pt}{20pt} 0 & \mbox{ if $n$ is odd;} \end{array} \right.\\
(\zeta^{-1})_-(M_\alpha)  & = \begin{cases}
{\displaystyle (-1)^{k(\alpha)+\ipartfrac{k_o(\alpha)}}  \binom{-1/2}{ \ipart{k_o(\alpha)}} } & \text{ if $a_1$ is odd,} \\
\rule{0pt}{20pt}  0   & \text{ if $a_1$ is even;}
\end{cases}\\
\rule{0pt}{20pt}(\zeta^{-1})_+(M_\alpha)  & = \begin{cases}
{\displaystyle (-1)^{k(\alpha)+\frac{k_o(\alpha)}{2}}
 \binom{-1/2}{ k_o(\alpha)/2} } & \text{ if $n$ is even,} \\
\rule{0pt}{20pt}  0       & \text{ if $n$ is odd;}
\end{cases}\\
\rule{0pt}{20pt}(\zeta^{-1})_-(F_\alpha) & = 
(-1)^{\ipartfrac{n+1}}\binom{\peakint(\tilde{\alpha})-1/2}{\ipartn}\,;\\
(\zeta^{-1})_+(F_\alpha) & = 
\left\{ \begin{array}{ll}
\rule{0pt}{20pt} {\displaystyle (-1)^{\frac{n}{2}}\binom{\peakaug(\alpha')-1/2}{n/2}}
& \text{if $n$ is even,} \\
\rule{0pt}{20pt} 0
& \text{if $n$ is odd.} 
\end{array}\right.
\end{align}
We mention that the convolution powers of the universal character are 
\begin{equation}
\zeta^m(M_\alpha)  =\binom{m}{k(\alpha)}\,,
\end{equation}
\begin{equation}
\zeta^m(F_\alpha)  =\binom{m+n-k(\alpha)}{n}\,,
\end{equation}
for any integer $m$~\cite[Formula (4.4)]{ABS}.



\end{document}